\documentclass{article}

\usepackage{color,graphicx,times}

\makeatother

\let \ttorg \tt \def \tt{\ttorg \obeyspaces}

\begin{document}

\date{}

\title{\Large\bf An Extended Bracket Polynomial for Virtual Knots and Links}
\author{Louis H. Kauffman\\ Department of Mathematics, Statistics \\ and Computer Science (m/c
249)    \\ 851 South Morgan Street   \\ University of Illinois at Chicago\\
Chicago, Illinois 60607-7045\\ $<$kauffman@uic.edu$>$\\}

\maketitle

\thispagestyle{empty}

\subsection*{\centering Abstract}

{\em This paper defines a new
invariant of virtual knots and flat virtual knots.
We study this invariant in two forms: the {\it extended bracket invariant} and the {\it arrow polyomial.}
The extended bracket polynomial takes the form of a sum of virtual graphs with polynomial coefficients. The arrow polynomial
is a polynomial with a finite number of variables for any given virtual knot or link. We show how the extended bracket polynomial can
be used to detect non-classicality and to estimate virtual crossing number and genus for virtual knots and links.}
\bigbreak

\noindent {{\bf Keywords}: virtual knot theory, Jones polynomial, bracket state sum, extended bracket invariant, arrow polynomial,
 virtual crossing number, virtual genus}

\noindent {{\bf Mathematics Subject Classification 2000}: 57M25, 57M27}

\section{Introduction}
Virtual knot theory is an extension of classical knot theory to stabilized embeddings of circles into thickened orientable surfaces of
genus possibly greater than zero. Classical knot theory is the case of genus zero. There is a diagrammatic theory for studying virtual knots and links, 
and this diagrammatic theory lends itself to the construction of numerous new invariants of virtual knots as well as extensions of known
invariants. 
\bigbreak

This paper defines a new invariant of virtual links and flat virtual links that we call the {\it extended bracket invariant.}
For a given virtual link $K$, the extended bracket invariant is denoted by $<<K>>$ and takes values in the module generated by 
isotopy classes of virtual $4$-regular graphs over the ring of Laurent polynomials $Q[A,A^{-1}]$ where $Q$ denotes the rational numbers. 
A {\it virtual graph} is represented in the plane via a choice of cyclic orders at its nodes. The virtual crossings in a virtual graph
are artifacts of the choice of placement in the plane, and we allow detour moves for consecutive sequences of virtual crossings just as in the
virtual knot theory. Two virtual graphs are {\it isotopic} if there is a combination of planar graph isotopies and detour moves that connect them.
See section 7 for the theory and calculations for this extended bracket polynomial. The extended bracket is defined by a state summation
with a new reduction relation on the states of the original bracket state sum. Virtual graphs are
defined in Section 4.1. 
\bigbreak

One of the applications of these new invariants is to the category of flat virtual knots and links. In Section 5 we review the definitions of flat
virtuals and we prove that {\it isotopy classes of long flat virtual knots embed in the isotopy classes of 
all long virtual knots} via the ascending map $A: LVF \longrightarrow LVK$. See Section 5 for this result. This embedding of the long flat virtual knots means
that there are many invariants of them  obtained by appling any invariant $Inv$ of virtual knots via the composition with the ascending map $A.$ This
situation is in direct contrast to closed long virtual flats, where it is not so easy to define invariants. The extended bracket polynonmial is an
invariant of both long virtual flats and closed virtual flats.
\bigbreak  

This paper is relatively self-contained, with Section 2 reviewing the definitions of virtual knots and links, Section 3 reviewing the surface interpretations 
of virtuals and Section 4 reviewing the definitions and properties of flat virtuals. Section 6 is a review of the bracket polynomial and the Jones polynomial
for virtual knots and links. We recall the virtualization construction that produces infinitely many non-trivial virtual knots with unit Jones polynomial.
These examples are of interest since there are no known examples of classical non-trivial knots with unit Jones polynomial. One may conjecture that
all non-trivial examples produced by virtualization are non-classical. It may be that the virtualization of some non-trivial classical knot is 
isotopic to a classical knot, but we have no evidence that this can happen. The extended bracket invariant developed in this paper may help in
deciding this question of non-classicality. 
\bigbreak

Section 5 proves that long flat virtual knots (virtual strings in Turaev's terminology) embed into long virtual knots.
This result is a very strong tool for discriminating long virtual knots from one another. We give examples and computations using this technique
by applying invariants of virtual knots to corresponding ascending diagram images of flat virtual knots. We also recall the {\it odd writhe} of a virtual
knot -- an invariant that is a kind of self-linking number for virtual knots. A long virtual can be closed just as a long knot can be closed. It is 
a remarkable property of long virtuals (both flat and not flat) that there are non-trivial long examples whose closures are trivial. Thus one would 
like to understand the kernel of the closure mapping from long virtuals to virtuals both in the flat and non-flat categories. It is hoped that the 
invariants discussed in this paper will further this question.
\bigbreak

In Section 7 we give the definition of the extended bracket invariant $<<K>>$ and a number of examples of its calculation. These examples 
include a verification of the non-classicality of the simplest example of virtualization,
a verification that the Kishino diagram \cite{KIS} and the flat Kishino diagram are non-trivial, and a verification that a particular flat diagram is non-trivial.
We use the extended bracket state sum to prove that an infinite family of single crossing virtualizations of classical diagrams are non-trivial and 
non-classical. The extended bracket is an invariant of flat diagrams by taking the specialization of its parameters so that $A=1$ and $\delta = -2.$
See Section 7 for the details. We give an example showing that the extended bracket can detect a long virtual knot whose closure is trivial. This is a 
capability that is beyond the reach of the Jones polynomial.
\bigbreak 

Section 8 proves two estimates for the virtual crossing number $VC(K)$ for a virtual link $K.$ The virtual crossing number is the least number of virtual 
crossings in any planar diagram that represents the virtual link. It is an isotopy invariant of the virtual link. Our estimates are based on the virtual crossing
numbers of the graphs that appear in the extended bracket state sum. The virtual crossing numbers of these graphs are usually easier to determine since the 
graphs have a more rigid structure than the links. We combine this observation with the fact that an adequate or semi-adequate link (one whose $A$ or $A^{-1}$ states do
not have any self-touching sites) has highest or lowest degree terms that can be pinpointed without calculating the entire state sum. This means that one can, for
such links, estimate the virtual crossing number without calculating the entire extended bracket state sum. We use this method to calculate virtual crossing 
numbers for the examples in Section 7 and we give an infinite collection of virtual links $L(n)$ whose virtual crossing number is $n = 1,2,3,\cdots$ and such that
each $L(n)$ can be represented by an embedding in a thickened torus. Finally, we give an infinite collection of virtual knots $K(n,m)$ with minimal genus $n+1$ and virtual
crossing number $n+m.$
\bigbreak 

Section 9 constructs a {\it simple extended bracket invariant}  denoted by ${\cal A}[K].$ We also refer to ${\cal A}[K]$  as the {\it arrow polynomial}
\cite{DyeKauff}. The arrow polynomial has infinitely many variables and integer coefficients. This invariant of regular isotopy of virtual knots and links
is obtained by the same method as the extended bracket, but we weaken the state reductions so that the reverse oriented smoothings each become two individual
graphical vertices. These vertices do not all disappear in the reduction process and one obtains reduced states that are disjoint unions 
of decorated circle graphs. These non-trivial circle graphs are denoted by commuting algebraic variables $K_{n}$ so that ${\cal A}[K]$ is a polynomial 
in the variables $A$, $A^{-1}$ and the $K_{n}.$ The invariant ${\cal A}[K]$ is quite strong and can be determined by a computer program that we describe
in this section. Calculations in this section are compared with the graphical calculations of the earlier part of the paper. The arrow polynomial has been
discovered independently by Miyazawa, via a different definition in \cite{Miyazawa2}. In fact, the work in the present paper was inspired by earlier work
of Miyazawa and Naoko Kamada. The papers \cite{Miyazawa2,KamadaMiya}, connect directly with the present work. The relationships between
the extended bracket invariant, the arrow polynomial and Miyazawa's work will be the subject of further investigation. 
In Section 10 we discuss how the arrow polynomial can be used in conjunction with the extended bracket to compute estimates on the genus of a virtual link or flat virtual link.
We give specific examples of genus computations at the end of the section.
Section 11 is a short discussion of open problems and future directions.
\bigbreak

\section{Virtual Knot Theory}
Knot theory
studies the embeddings of curves in three-dimensional space.  Virtual knot theory studies the  embeddings of curves in thickened surfaces of arbitrary
genus, up to the addition and removal of empty handles from the surface. Virtual knots have a special diagrammatic theory, described below,
that makes handling them
very similar to the handling of classical knot diagrams. Many structures in classical knot
theory generalize to the virtual domain.
\bigbreak  

In the diagrammatic theory of virtual knots one adds 
a {\em virtual crossing} (see Figure 1) that is neither an over-crossing
nor an under-crossing.  A virtual crossing is represented by two crossing segments with a small circle
placed around the crossing point. 
\bigbreak

Moves on virtual diagrams generalize the Reidemeister moves for classical knot and link
diagrams.  See Figure 1.  One can summarize the moves on virtual diagrams by saying that the classical crossings interact with
one another according to the usual Reidemeister moves while virtual crossings are artifacts of the attempt to draw the virtual structure in the plane. 
A segment of diagram consisting of a sequence of consecutive virtual crossings can be excised and a new connection made between the resulting
free ends. If the new connecting segment intersects the remaining diagram (transversally) then each new intersection is taken to be virtual.
Such an excision and reconnection is called a {\it detour move}.
Adding the global detour move to the Reidemeister moves completes the description of moves on virtual diagrams. In Figure 1 we illustrate a set of local
moves involving virtual crossings. The global detour move is
a consequence of  moves (B) and (C) in Figure 1. The detour move is illustrated in Figure 2.  Virtual knot and link diagrams that can be connected by a finite 
sequence of these moves are said to be {\it equivalent} or {\it virtually isotopic}.
\bigbreak

\begin{figure}[htb]
     \begin{center}
     \begin{tabular}{c}
     \includegraphics[width=10cm]{F1.EPSF}
     \end{tabular}
     \caption{\bf Moves}
     \label{Figure 1}
\end{center}
\end{figure}

\begin{figure}[htb]
     \begin{center}
     \begin{tabular}{c}
     \includegraphics[width=10cm]{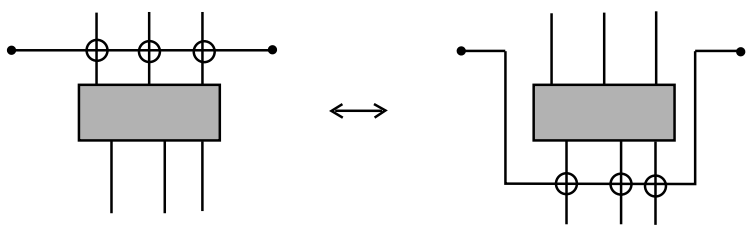}
     \end{tabular}
     \caption{\bf Detour Move}
     \label{Figure 2}
\end{center}
\end{figure}

\begin{figure}[htb]
     \begin{center}
     \begin{tabular}{c}
     \includegraphics[width=10cm]{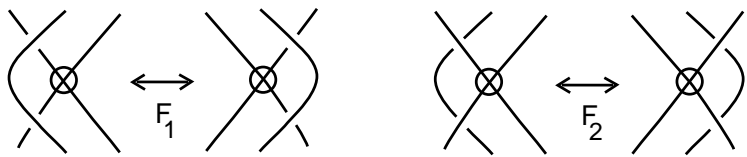}
     \end{tabular}
     \caption{\bf Forbidden Moves}
     \label{Figure 3}
\end{center}
\end{figure}

Another way to understand virtual diagrams is to regard them as representatives for oriented Gauss codes \cite{GPV}, \cite{VKT,SVKT} 
(Gauss diagrams). Such codes do not always have planar realizations. An attempt to embed such a code in the plane
leads to the production of the virtual crossings. The detour move makes the particular choice of virtual crossings 
irrelevant. {\it Virtual isotopy is the same as the equivalence relation generated on the collection
of oriented Gauss codes by abstract Reidemeister moves on these codes.}  
\bigbreak

Figure $3$ illustrates the two {\it forbidden moves}. Neither of these follows from Reidmeister moves plus detour move, and 
indeed it is not hard to construct examples of virtual knots that are non-trivial, but will become unknotted on the application of 
one or both of the forbidden moves. The forbidden moves change the structure of the Gauss code and, if desired, must be 
considered separately from the virtual knot theory proper. 
\bigbreak

\section{Interpretation of Virtuals Links as Stable Classes of Links in  Thickened Surfaces}
There is a useful topological interpretation \cite{VKT,DVK} for this virtual theory in terms of embeddings of links
in thickened surfaces.  Regard each 
virtual crossing as a shorthand for a detour of one of the arcs in the crossing through a 1-handle
that has been attached to the 2-sphere of the original diagram.  
By interpreting each virtual crossing in this way, we
obtain an embedding of a collection of circles into a thickened surface  $S_{g} \times R$ where $g$ is the 
number of virtual crossings in the original diagram $L$, $S_{g}$ is a compact oriented surface of genus $g$
and $R$ denotes the real line.  We say that two such surface embeddings are
{\em stably equivalent} if one can be obtained from another by isotopy in the thickened surfaces, 
homeomorphisms of the surfaces and the addition or subtraction of empty handles (i.e. the knot does not go through the handle).

\begin{figure}
     \begin{center}
     \begin{tabular}{c}
     \includegraphics[width=10cm]{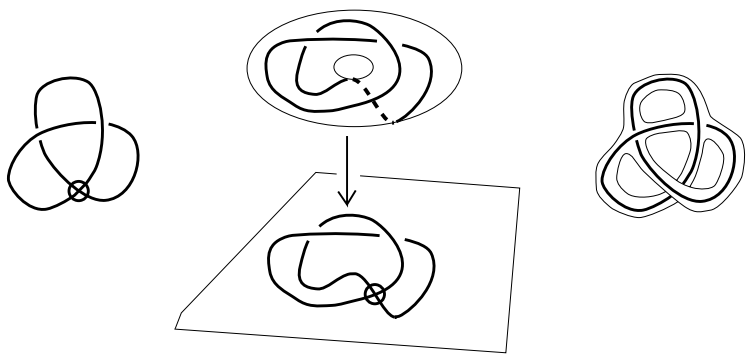}
     \end{tabular}
     \caption{\bf Surfaces and Virtuals}
     \label{Figure 4}
\end{center}
\end{figure}

\noindent We have the
\smallbreak
\noindent
{\bf Theorem 1 \cite{VKT,DKT,DVK,Carter}.} {\em Two virtual link diagrams are isotopic if and only if their corresponding 
surface embeddings are stably equivalent.}  
\smallbreak
\noindent
\bigbreak  

\noindent In Figure 4 we illustrate some points about this association of virtual diagrams and knot and link diagrams on surfaces.
Note the projection of the knot diagram on the torus to a diagram in the plane (in the center of the figure) has a virtual crossing in the 
planar diagram where two arcs that do not form a crossing in the thickened surface project to the same point in the plane. In this way, virtual 
crossings can be regarded as artifacts of projection. The same figure shows a virtual diagram on the left and an ``abstract knot diagram" \cite{Kamada3,Carter} on the right.
The abstract knot diagram is a realization of the knot on the left in a thickened surface with boundary and it is obtained by making a neighborhood of the 
virtual diagram that resolves the virtual crossing into arcs that travel on separate bands. The virtual crossing appears as an artifact of the
projection of this surface to the plane. The reader will find more information about this correspondence \cite{VKT,DKT} in other papers by the author and in
the literature of virtual knot theory.
\bigbreak
 
\section{Flat Virtual Knots and Links}
Every classical knot or link diagram can be regarded as a $4$-regular plane graph with extra structure at the 
nodes. This extra structure is usually indicated by the over and under crossing conventions that give
instructions for constructing an embedding of the link in three dimensional space from the diagram.  If we take the flat diagram
without this extra structure then the diagram is the shadow of some link in three dimensional space, but the weaving of that link is not 
specified. It is well known that if one is allowed to apply the Reidemeister moves to such a shadow (without regard to the types
of crossing since they are not specified) then the shadow can be reduced to a disjoint union of circles. This reduction is 
no longer true for virtual links. More precisely, let a {\em flat virtual diagram} be a diagram with {\it virtual crossings} as we have
described them and {\em flat crossings} consisting in undecorated nodes of the $4$-regular plane graph. Two flat virtual diagrams are {\em equivalent} if
there is a  sequence of generalized flat Reidemeister moves (as illustrated in Figure 1) taking one to the other. A generalized
flat Reidemeister move is any move as shown in Figure 1 where one ignores the over or under crossing structure.
Note that in studying flat virtuals the rules for changing virtual crossings among themselves and the rules for changing
flat crossings among themselves are identical. Detour moves as in Figure 1C are available for virtual crossings
with respect to flat crossings and {\it not} the other way around. The analogs of the forbidden moves of Figure 3 remain forbidden 
when the classical crossings are replaced by flat crossings. 
\bigbreak

The theory of flat virtual knots and links is identical to the theory of all oriented Gauss codes (without over or under information)
modulo the flat Reidemeister moves. Virtual crossings are an artifact of the realization of the flat diagram in the plane.
In Turaev's work \cite{VST} flat virtual knots and links are called {\it virtual strings}. See also recent papers of Roger Fenn \cite{Fenn1,Fenn2} for other
points of view about flat virtual knots and links.
\bigbreak

We shall say that a virtual diagram {\em overlies} a flat diagram if the virtual diagram is obtained from the flat diagram by
choosing a crossing type for each flat crossing in the virtual diagram. To each virtual diagram $K$ there is an associated 
flat diagram $F(K)$, obtained by forgetting the extra structure at the classical crossings in $K.$ Note that if $K$ and $K'$
are isotopic as virtual diagrams, then $F(K)$ and $F(K')$ are isotopic as flat virtual diagrams. Thus, if we can
show that $F(K)$ is not reducible to a disjoint union of circles, then it will follow that $K$ is a non-trivial virtual link.
The flat virtual diagrams present a challenge for the
construction of new invariants. They are fundamental to the study of virtual knots.  A virtual knot is necessarily non-trivial if its flat projection
is a non-trivial flat virtual diagram.  We wish to be able to determine when a given virtual  link is isotopic to a classcal
link. The reducibility or irreducibility of the underlying flat diagram is the first  obstruction to such an equivalence.
\bigbreak 

\noindent {\bf Definition.}
A {\it virtual graph} is a flat virtual diagram where the classical flat crossings are not subjected to the flat Reidemeister moves.
Thus a virtual graph is a $4-regular$ graph that is represented in the plane via a choice of cyclic orders at its nodes. The virtual crossings
are artifacts of this choice of placement in the plane, and we allow detour moves for consecutive sequences of virtual crossings just as in the
virtual knot theory. Two virtual graphs are {\it isotopic} if there is a combination of planar graph isotopies and detour moves that connect them.
The theory of virtual graphs is equivalent to the theory of $4$-regular graphs on oriented surfaces, taken up to empty handle stabilization, in 
direct analogy to our description of the theory of virtual links and flat virtual links.
\bigbreak

\section{Long Knots and Long Flats}
A long knot or link is a $1-1$ tangle. It is a tangle with one input end and one output end. In between one has, in the diagram, 
any knotting or linking, virtual or classical. Classical long knots (one component) carry essentially the same topological information
as their closures. In particular, a classical long knot is knotted if and only if its closure is knotted. This statement is false for virtual
knots. An example of the phenomenon is shown in Figure 5.
\bigbreak

The long knots $L$ and $L'$ shown in Figure 5 are non-trivial in the virtual category. Their closures, obtained by attaching the ends together are
unknotted virtuals. Concomittantly, there can be a multiplicity of long knots associated to a given virtual knot diagram, obtained by cutting an
arc from the diagram and creating a $1-1$ tangle. {\it It is a fundamental problem to determine the kernel of the closure mapping from long
virtual knots to virtual knots.} In Figure 5, the long knots $L$ and $L'$ are trivial as welded long knots (where the first forbidden move is allowed). The
obstruction to untying them as virtual long knots comes from the first forbidden move. The matter of proving that $L$ and $L'$ are non-trivial distinct long
knots is difficult by direct attack. There is a fundamental relationship betweem long flat virtual knots and long virtual knots that can be used to
see it.
\bigbreak

Let $LFK$ denote the set of long flat virtual knots and let $LVK$ denote the set of long virtual knots.
We define $$A:LFK \longrightarrow LVK$$ by letting $A(S)$ be the ascending long virtual knot diagram associated with 
the long flat virtual diagram $S.$ That is, $A(S)$ is obtained from $S$ by traversing $S$ from its left end to its right end and creating
a crossing at each flat crossing so that {\it one passes under each crossing before passing over that crossing}. Virtual crossings are not changed by this
construction. The idea of using the ascending diagram to define invariants of long flat virtuals is exploited in \cite{SW}. The following result is new.
\bigbreak

\noindent {\bf Long Flat Embedding Theorem 2.} The mapping $A:LFK \longrightarrow LVK$ is well-defined on the corresponding isotopy classes of diagrams, 
and it is injective. Hence {\it long flat virtual knots embed in the class of long virtual knots.}
\bigbreak

\noindent {\bf Proof.} It is easy to see that if two flat long diagrams $S,T$ are virtually isotopic then $A(S)$ and $A(T)$ are isotopic 
long virtual knots. We leave the verification to the reader. To see the injectivity, define 
$$Flat:LVK \longrightarrow LFK$$ by letting $Flat(K)$ be the long flat diagram obtained from $K$ by flattening all the classical crossings in $K.$
By definition $Flat(A(S)) = S.$ Note that the map $Flat$ takes isotopic long virtual knots to isotopic flat virtual knots since
any move using the under and over crossings is a legitimate move when this distinction is forgotten. This proves injectivity and
completes the proof of the Theorem. //
\bigbreak

The long flat theorem is easy to prove, and it has many good consequences. First of all, let $Inv(K)$ denote any invariant
of long virtual knots. ($Inv(K)$ can denote a polynomial invariant, a number, a group, quandle, biquandle or other invariant structure.) Then we can define
$$Inv(S)$$ for any long flat knot $S$ by the formula $$Inv(S) = Inv(A(S)),$$ and this definition yields an invariant of long flat knots.
\bigbreak

\begin{figure}
     \begin{center}
     \begin{tabular}{c}
     \includegraphics[width=8cm]{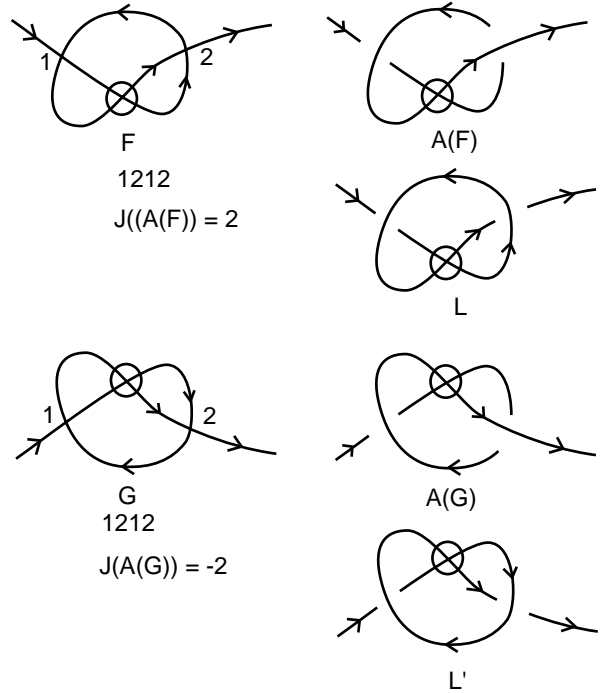}
     \end{tabular}
     \caption{\bf Ascending Map}
     \label{Figure 5}
\end{center}
\end{figure}

\begin{figure}
     \begin{center}
     \begin{tabular}{c}
     \includegraphics[width=8cm]{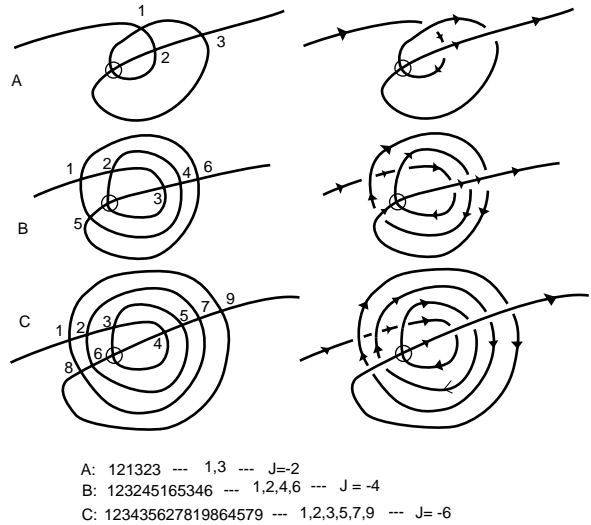}
     \end{tabular}
     \caption{\bf Spiral Examples}
     \label{Figure 6}
\end{center}
\end{figure}

\noindent {\bf Using the Odd Writhe J(K).} View Figure 5. We show the long flat $F$ and its image under the
ascending map, $A(F)$ are non-trivial. In fact, $A(F)$ is non-trivial and non-classical. One computes that $J(A(F))$ is
non-zero where $J(K)$ denotes the {\it odd writhe} of $K.$ The odd writhe \cite{SL} is the sum of the signs of the odd crossings. A crossing is {\it odd} if 
it flanks an odd number of symbols in the Gauss code of the diagram. Classical diagrams have zero odd writhe. 
In this case, the flat Gauss code for $F$ is $1212$ with both  crossings odd. Thus we see from the figure that $J(A(F)) = 2.$
Thus $A(F)$ is non-trivial, non-classical and inequivalent to its mirror image.
Once we check that $A(F)$ is non-trivial, we know that the flat knot $F$ is non-trivial, and from this we conclude that the long virtual
$L$ is also non-trivial. Note that $J(L) = 0,$ so we cannot draw this last conclusion directly from $J(L).$ This same figure illustrates a long flat $G$
that is obtained by reflecting $F$ in a horizontal line. Then, as the reader can calculate from this figure, $J(A(G)) = -2.$ Thus $F$ and $G$ 
are distinct non-trivial long flats. We conclude from these arguments that the long virtual knots $L$ and $L'$ in Figure 5 are both non-trivial, and that
$L$ is not virtually isotopic to $L'$ (since such an isotopy would give an isotopy of $F$ with $G$ by the flattening map).
\bigbreak

A second class of examples is shown in Figure 6. Here the examples labeled $A$, $B$ and $C$ are part of an infinite family of long flat virtuals
whose ascending long virtuals have odd writhe $J = -2n$ for $n = 1, 2, 3, \cdots$. This gives an infinite family of distinct long flat virtuals such that each
one has trivial closure. In Figure 6 we give the flat Gauss code for each example, then list the odd crossings and the value of the invariant $J.$
\bigbreak

\section{Review of the Bracket Polynomial for Virtual Knots}

In this section we recall how the bracket state summation model \cite{K} for the Jones polynomial \cite{Jones,WITT} is defined for virtual knots
and links.  In the next section we give an extension of this model using orientation structures on the states of the bracket expansion.
The extension is also an invariant of flat virtual links.
\bigbreak

We call a diagram in the plane 
{\em purely virtual} if the only crossings in the diagram are virtual crossings. Each purely virtual diagram is equivalent by the
virtual moves to a disjoint collection of circles in the plane.
\bigbreak

A state $S$ of a link diagram $K$ is obtained by
choosing a smoothing for each crossing in the diagram and labelling that smoothing with either $A$ or $A^{-1}$
according to the convention that a counterclockwise rotation of the overcrossing line sweeps two 
regions labelled $A$, and that a smoothing that connects the $A$ regions is labelled by the letter $A$. Then, given
a state $S$, one has the evaluation $<K|S>$ equal to the product of the labels at the smoothings, and one has the 
evaluation $||S||$ equal to the number of loops in the state (the smoothings produce purely virtual diagrams).  One then has
the formula
$$<K> = \Sigma_{S}<K|S>d^{||S||-1}$$
where the summation runs over the states $S$ of the diagram $K$, and $d = -A^{2} - A^{-2}.$
This state summation is invariant under all classical and virtual moves except the first Reidemeister move.
The bracket polynomial is normalized to an
invariant $f_{K}(A)$ of all the moves by the formula  $f_{K}(A) = (-A^{3})^{-w(K)}<K>$ where $w(K)$ is the
writhe of the (now) oriented diagram $K$. The writhe is the sum of the orientation signs ($\pm 1)$ of the 
crossings of the diagram. The Jones polynomial, $V_{K}(t)$ is given in terms of this model by the formula
$$V_{K}(t) = f_{K}(t^{-1/4}).$$
\noindent This definition is a direct generalization to the virtual category of the  
state sum model for the original Jones polynomial. It is straightforward to verify the invariances stated above.
In this way one has the Jones polynomial for virtual knots and links.
\bigbreak

\noindent We have \cite{DVK} the  
\smallbreak
\noindent
{\bf Theorem 3.} {\em To each non-trivial
classical knot diagram of one component $K$ there is a corresponding  non-trivial virtual knot diagram $Virt(K)$ with unit
Jones polynomial.} 
\bigbreak

\noindent {\bf Proof Sketch.} This Theorem is a key ingredient in the problems involving virtual knots. Here is a sketch of its proof.
The proof uses two invariants of classical knots and links that generalize to arbitrary virtual knots and links.
These invariants are the {\em Jones polynomial} and the {\em involutory quandle} denoted by the notation
$IQ(K)$ for a knot or link $K.$ 
\bigbreak

Given a
crossing $i$ in a link diagram, we define $s(i)$ to be the result of {\em switching} that crossing so that the undercrossing arc
becomes an overcrossing arc and vice versa. We define the {\em virtualization}
$v(i)$ of the crossing by the local replacement indicated in Figure 7. In this figure we illustrate how, in virtualization, 
the  original crossing is replaced by a crossing that is flanked by two virtual crossings. When we smooth the two
virtual crossings in the virtualization we obtain the original knot or link diagram with the crossing switched.
\bigbreak

Suppose that $K$ is a (virtual or classical) diagram with a classical crossing labeled $i.$  Let $K^{v(i)}$ be the diagram
obtained from $K$ by virtualizing the crossing $i$ while leaving the rest of the diagram just as before. Let $K^{s(i)}$ be
the diagram obtained from $K$ by switching the crossing $i$ while leaving the rest of the diagram just as before. Then it
follows directly from the expansion formula for the bracket polynomial that $$V_{K^{s(i)}}(t) = V_{K^{v(i)}}(t).$$ 
\noindent As far as the Jones
polynomial is concerned, switching a crossing and virtualizing a crossing look the same. We can start with a classical knot diagram $K$
and choose a subset $S$ of crossings such that the diagram is unknotted when these crossings are switched. Letting $Virt(K)$ denote the virtual knot
diagram obtained by virtualizing each crossing in $S$, it follows that the Jones polynomial of $Virt(K)$ is equal to unity, the Jones polynomial
of the unknot. Nevertheless, if the original knot $K$ is knotted, then the virtual knot $Virt(K)$ will be non-trivial. We outline the argument for this
fact below.
\bigbreak 
 
The involutory quandle \cite{KNOTS} is an algebraic invariant
equivalent to the fundamental group of the double branched cover of a knot or link in the classical case. In this algebraic
system one associates a generator of the algebra $IQ(K)$ to each arc of the diagram $K$ and there is a relation of the form
$c = ab$ at each crossing, where $ab$ denotes the (non-associative) algebra product of $a$ and $b$ in $IQ(K).$ See Figure 8.
In this figure we have illustrated the fact that  $$IQ(K^{v(i)}) = IQ(K).$$
As far as the involutory quandle is concerned, the original crossing and the virtualized crossing look the same.
\bigbreak

If a classical knot is actually knotted, then its involutory quandle is non-trivial \cite{W}. Hence if we start
with a non-trivial classical knot and virtualize any subset of its crossings we obtain a virtual knot that is still 
non-trivial. There is a subset $A$ of the crossings of a classical knot $K$ such that the knot $SK$ obtained by
switching these crossings is an unknot.  Let $Virt(K)$ denote the virtual diagram obtained from $A$ by virtualizing
the crossings in the subset $A.$  By the above discussion, the Jones polynomial of $Virt(K)$ is the same
as the Jones polynomial of $SK$, and this is $1$ since $SK$ is unknotted. On the other hand, the $IQ$ of $Virt(K)$ is the 
same as the $IQ$ of $K$, and hence if $K$ is knotted, then so is $Virt(K).$   We have shown that $Virt(K)$ is a non-trivial
virtual knot with unit Jones polynomial.  This completes the proof of the Theorem. //
\bigbreak

\begin{figure}
     \begin{center}
     \begin{tabular}{c}
     \includegraphics[width=4cm]{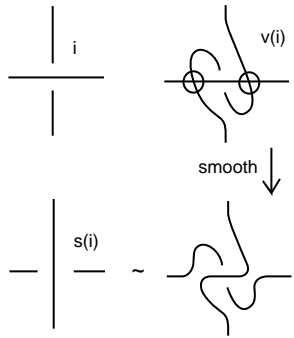}
     \end{tabular}
     \caption{\bf Switch and Virtualize}
     \label{Figure 7}
\end{center}
\end{figure}

\begin{figure}
     \begin{center}
     \begin{tabular}{c}
     \includegraphics[width=4cm]{F8.EPSF}
     \end{tabular}
     \caption{\bf IQ(Virt)}
     \label{Figure 8}
\end{center}
\end{figure}

See Figure 32 for an example of a virtualized trefoil, the simplest example of a non-trivial virtual knot with unit Jones polynomial.
More work is needed to prove that the virtual knot $T$ in Figure 32 is not classical. In the next section we will give a proof of this fact
by using an extension of the bracket polynomial.
\bigbreak

It is an open problem whether there are classical knots (actually knotted) having unit Jones polynomial. (There are linked links whose linkedness is unseen
\cite{MT,EKT} by the Jones polynomial.) If there exists a classical knot with unit Jones polynomial, then one of the knots $Virt(K)$ produced by this Theorem
may be isotopic to  a classical knot.  Such examples are guaranteed to be non-trivial, but they are usually also not classical. 
We do not know at this writing whether all such virtualizations of non-trivial classical knots, yielding virtual knots with unit Jones polynomial, are 
non classical. It is
an intricate task to verify that specific examples of 
$Virt(K)$ are not classical. This has led to an investigation of new invariants for virtual knots. In this way the search for classical knots with unit Jones
polynomial expands to exploration of the structure of the infinite collection of virtual knots with unit Jones polynomial.
\bigbreak

\begin{figure}
     \begin{center}
     \begin{tabular}{c}
     \includegraphics[width=3cm]{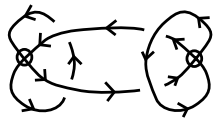}
     \end{tabular}
     \caption{\bf Kishino Diagram}
     \label{Figure 9}
\end{center}
\end{figure}

In Figure 9 we show the {\it Kishino diagram} $K.$
This diagram has unit Jones polynomial and its fundamental group is infinite cyclic. The Kishino diagram was discovered by Kishino in \cite{KIS}.
Many other invariants of virtual knots 
fail to detect the Kishino knot. Thus it has been a test case for examining new invariants. Heather Dye and the author \cite{MinSurf}
have used the bracket polynomial defined for knots and links in a thickened surface (the state curves
are taken as isotopy classes of curves in the surface) to prove the non-triviality and non-classicality of the Kishino diagram.
In fact, we have used this technique to show that knots with unit Jones polynomial obtained by a single virtualization are non-classical.
See the problem list by Fenn, Kauffman and Manturov \cite{VP} for other problems and proofs related to the Kishino diagram.
In the next section we describe a new extension of the bracket polynomial that can be used to discriminate the Kishino diagram, and, in fact,
shows that its corresponding flat virtual knot is non-trivial.
\bigbreak

\section{An Extended Bracket Polynomial for Virtual and Flat Virtual Knots and Links}
This section describes a new invariant for virtual knots and links,  and for flat virtual knots and links. The construction of the invariant begins with the
oriented state summation of the bracket polynomial. This means that each local smoothing is either an oriented smoothing or a {\it disoriented
smoothing} as illustrated in Figures 10 and 11.
\bigbreak

\begin{figure}
     \begin{center}
     \begin{tabular}{c}
     \includegraphics[width=6cm]{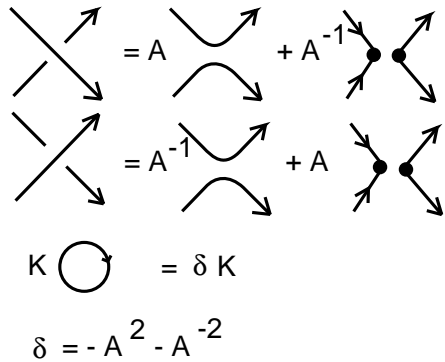}
     \end{tabular}
     \caption{\bf Oriented Bracket Expansion}
     \label{Figure 10}
\end{center}
\end{figure}

In Figure 10 we illustrate the oriented bracket expansion for both positive and negative crossings in the link diagram. An oriented crossing can
be smoothed in the oriented fashion or the disoriented fashion as shown in Figure 10. We refer to these smoothings as {\it oriented} and {\it disoriented}
smoothings. To each smoothing we make an associated configuration that will be part of the extended state summation. The configuration associated to the
oriented smoothing is that smoothing itself. The configuration associated to the disoriented smoothing is obtained by applying the reduction rules described below. 
The reduction process yields reduced states where the remaining disoriented smoothings are replaced by graphical vertices.  See Figures 11 and 14. The
extended bracket state summation is defined by the formula:
$$<<K>> = \Sigma_{S}<K|S>d^{||S||-1} [S]$$
where $S$ runs over the oriented bracket states of the diagram, $<K|S>$ is the usual product of vertex weights as in the 
standard bracket polynomial, and $[S]$ is a sum of reduced states (flat virtual
graphs) associated with the state
$S$ that is obtained by  rules that we describe below.
\bigbreak

The square brackets around $S$ in the state summation denote its replacement by a sum of {\it reduced states}, using the replacements of Figures
11 and 12.  The replacement procedure will 
be explained in detail below. We take each state and apply reducing rules to it, sometimes creating a multiplicity of graphs from the given state 
graph. The disoriented vertex in a state is regarded as a special 4-regular vertex at the end of the process of reduction. When the state is fully reduced we draw
the remaining disoriented sites as $4$-regular vertices in the reduced graph, or leave them as disoriented sites. At this stage, disoriented sites and
$4$-regular vertices are interchangeable. We  prove below that virtual diagrams related by Reidemeister moves have the
same set of reduced graphs. We then prove that the extended bracket polynomial is an invariant of regular isotopy for virtual links.
\bigbreak

The extended bracket of
$K$ takes values in the module generated by $4$-regular virtual graphs with coefficients in the Laurent polynomial ring
$Q[A,A^{-1}].$ It is an element in the module over $Q[A,A^{-1}]$ generated by the planar isotopy classes of these virtual graphs. 
\bigbreak

The summation of reduced states $[S]$ associated with a state
$S$ is obtained by following, in an order we describe below, the reduction rules of Figures 11 and 12.  Virtual crossings are not changed by the state reduction.  
\bigbreak

It is convenient to use, in one step, certain consequences of the basic replacements of Figure 11. In Figure 12 we show three {\it special replacements} 
labeled $A,B,C$. The special
replacements are consequences of the basic replacements of Figure 11. We will see that a formalized reduction of states via these replacements will
ensure the invariance of $<<K>>$ under the Reidemeister and detour moves. Before giving the precise reduction rules, we make a remark about how they
are derived. 
\bigbreak 

\noindent {\bf Remark.} In order to see the motivation for the reduction rules below, we recommend that the reader examine the state
expansions for the Reidemeister moves as illustrated in Figures 13, 23, 24 and 25.  It is necessary to allow certain reductions in order to obtain invariance
under the Reidemeister moves. By examining the expansions in Figures 13, 23, 24 and 25, the reader will see how the definition of the reductions originates in
tracking the combinatorics of these moves.
It should be clear (and will be discussed below) that move 2 of Figure 11 is needed for the first Reidemeister move, and that the special moves A, B and C
of Figure 12 are
needed for the second and third Reidemeister moves. In particular, moves A and B are needed for the second Reidemeister move, and move C is needed for the third
Reidemeister move. {\it Figure 26 shows, that the moves 2, B and C will be incompatible unless we also have available the move 3 of Figure 11.} Then, 
adding move 3, we see easily that the moves A, B and C are consequences of the rules in Figure 11. Nevertheless, we find that move C must be given
precedence, as explained below.
This is the basic analysis that leads to our systematization of the state reductions in terms of the moves of Figures 11, 12 and the rules given below.
\bigbreak

\noindent {\bf The rules for forming replacements on
a state of the extended bracket polynomial:}

\begin{enumerate}
\item First consider all the oriented smoothings in the state. Each oriented smoothing site in $S$ is left as an oriented smoothing in $[S]$.
\item Perform all single oriented loop simplifications, labeled type $1$ in Figure 11. In this replacement, a disjoint  oriented loop in a state does
not contribute to the reduced state, and is removed.
\item Search the state for all replacements of type $C$ of Figure 12, and make these replacements. Repeat the previous step if necessary.
\item Search the state for all replacments of type $2$ and type $3$  of Figure 11, and make these replacements. Repeat the previous steps if necessary --
that is, when oriented loops or type $C$ configurations appear, repeat steps 2 and 3.
If there are a multiplicity of replacements of type $3$ available, as shown generically in Figure 14, make each replacement
and take the sum of all of them divided by the number of such replacements.
\item When there are no more simplifications available, replace each remaining disoriented smoothing with a graphical vertex as shown in second item from 
the top of Figure 11.
\end{enumerate}
\bigbreak

\noindent {\bf Remark.} The last replacement of disoriented smoothings (paired cusps) with graphical vertices is a notational device to indicate 
that the state is in reduced form. Sometimes it is convenient to simply note that a collection of states is reduced, and not replace the paired cusps with
graphical vertices (e.g. see Figure 41 where both methods of displaying reduced states are used). In all cases, we regard each paired cusp as equivalent to a 4-valent vertex in the plane, and each state is treated
(reduced or unreduced) as a virtual graph, taken up to planar isotopy and virtual detour moves. See Figures 30 and 31 where we illustrate the effect of virtual
detour moves on the representation of states of two diagrams. In these figures the detour moves are indicated by a tilde between the corresponding state 
diagrams.
\bigbreak

The reduction rules do not explicitly mention the special replacements of type $A$ and $B$ from Figure 12.  Nevertheless, we need, for the sake of invariance under the Reidemeister moves, to analyze states that are related by reductions of all three types
$A$, $B$ and $C.$  We now explain how these replacements are related to invariance under the Reidemeister moves. 
View Figure 13 and note that we obtain framing behaviour under the first Reidemeister
move just as in the case of the classical bracket polynomial. View Figure 12 and note that each of the special replacments shown is 
a consquence of the reduction rules of Figure 11. In the reduction rules we have made the special replacement $C$ of Figure 12 (which is locally 
recognizable in a diagram) take precedence over the application of the simpler rules of Figure 11. The reason for this precedence choice is to eliminate
other reduction possibilities that can occur if this ordering is not followed. In Figure 21 we illustrate how a reduction with different cusp pairings 
can occur if $C$ is not
performed first. 
\bigbreak

\noindent {\bf Definition.} Given a state $S$ of an oriented virtual diagram $K,$ let $Red(S)$ denote the set of all reduced states obtained from $S$ by the
reduction process described above. (We have already remarked that there may be a multiplicity of distinct reduced states associated with a given state $S.$)
\bigbreak

The special replacements $A$ and $B$ have the property that they do not, by themselves, generate a multiplicity of reductions. We have the following lemma.
\bigbreak

\noindent {\bf Reduction Lemma 1.} Let a state $S$ of an oriented virtual diagram $K$ have the initial configuration for $A$ or for $B$ in Figure 12,
and let $S'$ be the state obtained from
$S$ by applying the reduction of type $A$ or $B.$ Then (see the definition above) $$Red(S) = Red(S').$$
\smallbreak

\noindent {\bf Proof.} We must consider how a diagram with a type $A$ or type $B$ configuration will be changed by the reduction rules.
View Figures 15, 16 and 17. These figures illustrate a generic cases involving the type $A$ replacement of Figure 12. Each arrow in Figure 15 is an
application of the basic replacement of type 3 in Figure 11. Note that we indicate type replacements of type 3 either by a label that refers to
an edge in the diagram, or by the label $[3],$ referring to type 3. While, in Figure 15, there are distinct intermediate reductions, 
there is a unique final reduction described
by one application of the type $A$ replacement of Figure 12. This demonstrates that applying the type $A$ replacement will result in the same set of reduced states
as the intial diagram. Figure 16 illustrates the same property in the presence of a multiplicity of type A moves. Figure 17 shows a corresponding situation
where the state reduction can factor through types $[1]$ and A, or through types $[3], [1]$ and B.
Figures 18, 19 and 20  show the same properties for the type $B$ replacement. //
\bigbreak 

\noindent {\bf Reduction Lemma 2.} Let a state $S$ of an oriented virtual diagram $K$ have an intial configuration of type $C$ as in Figure 12.
Let $S'$ be the state obtained
by applying the reduction of type $C$ to this local configuration in $S.$ Then $$Red(S) = Red(S').$$
\bigbreak

\noindent {\bf Proof.} Without loss of generality, we can assume that the state $S$ has no simple oriented disjoint loops, since these will be eliminated by the
first reduction rule. The second reduction rule demands that type $C$ reduction be applied wherever it is applicable. We claim
that the reduction of $S$ by only type $C$ moves is unique and independent of the order of the application of these moves. To this purpose, view Figure 22.
In this figure we illustrate the generic possibilities that correspond to direct interactions of type $C$ reductions. {\it We note that the end results are 
independent of order of application of the type $C$ moves.} To see this, examine the vertices labeled $a$ and $b$ on the template for the $C$ move
in Figure 22. {\it It is not possible for a vertex to be
both type $a$ and type $b.$} This statement follows from the fact that the upper cusp in the type $b$ vertex is connected to 
the lower cusp through two cusps of type $a$ vertices. The upper cusp of a type $a$ vertex cannot connect to its lower cusp through only two other
vertices. If a $C$ move configuration shares vertices with another $C$ move configuration, then it can share a type $b$ vertex
with a vertex of type $b$ in the other configuration, or it can share a type $a$ vertex with another type $a$ vertex. 
In each case the end results of the $C$ reductions are independent of the order
of application of the $C$ moves. Cases of this sharing are illustrated in Figure 22. Since we know that the inital set of type $C$ reductions gives the 
same set of reductions independent of the order of operations, it follows that $Red(S) = Red(S').$ This completes the proof of the Lemma. //
\bigbreak

Now we analyze the Reidemeister moves. We precede the proof with an informal discussion.  We want to see that the sum of states after a Reidemeister move is the same as the sum before the Reidemeister move
is performed. To see the issue involved, view Figure 25 and particularly the three state-representatives in the rectangular box at the top of the figure.
Each small diagram represents all the states in the state sum that have the given local configuration. Each of these states will undergo a reduction process.
Concentrate on the middle diagram in the rectangle. This diagram has an instance of the reduction of type $C.$ Since we have given precedence to the type 
$C$ reduction, all the states associated with the middle diagram will reduce in the same way, and we can assume that the type $C$ reduction has been performed
at this diagram. The first diagram in the rectangle has a loop that will be eliminated by the special reduction of type $A$ (the loop is replaced 
by a factor of $d = -A^2 - A^{-2}$ in the state sum and the two cusps at the top and bottom of the diagram become neighbors). 
Each of the three diagrams is associated with the same set of reduced states. Thus the equation $d = -A^2 - A^{-2}$ assures us that these states all cancel from the sum of 
reduced states. 
\bigbreak

\noindent {\bf Theorem 4.} The extended bracket has framing behaviour under the first Reidemeister moves as shown in Figure 13 and it is invariant under
the second and third Reidemeister moves and under virtual detour moves. Hence the extended bracket state sum is a regular isotopy invariant of virtual
knots and links. 
\smallbreak

\noindent {\bf Proof.} We have already discussed the first Reidemeister move. Invariance under virtual detours is implicit in the definition of the 
state sum. In Figures 23 and 24 we show how the special removal ($B$ of Figure 12)  of a disoriented loop gives rise to invariance under the directly oriented
second Reidemeister move. In the state expansion, the disoriented loop is removed but multiplies its term by $d = -A^{2} - A^{-2}.$  The other two disoriented local configurations receive vertex weights
of $A^2$ and $A^{-2}$, and {\it each of these configurations has the same set of reduced states}. 
Thus these three configurations cancel each other from the state sum. This leaves the remaining local state with parallel arcs, 
and gives invariance under the directly oriented second Reidemeister move. Invariance under the  reverse oriented second Reidemeister move is shown in Figure 24. 
\bigbreak

Finally, view Figure 25 to see the structure of the third Reidemeister move.  In this figure we see
that  the special replacements ($B$ and $C$ of Figure 12) come into play to effect a cancellation of all terms shown in boxed format. The three terms on the
top row each have the same set of reduced states, and they add up to zero in the evaluation. 
The remaining terms in the sum are symmetrical with
respect to reflection in a vertical line and thus will contribute in the same way to the other side of the third Redeimeister move. This completes the
proof that  the extended bracket state sum is an invariant of regular isotopy. //
\bigbreak

\noindent {\bf Remark on State Reduction.} It is natural to ask for the combinatorial structure of the individual loops in a reduced state.
Note that we can regard a reduced state as decomposed into loops, or as a virtual graph obtained by binding the paired cusps of the remaining disoriented
vertices. Viewing the reduced state as a collection of loops, we can characterize the structure of each loop, independent of the pairings of the cusps. 
View Figures 48, 49 and 50. While these figures are designed to explain (for Section 9) a simpler reduction of states where we do not insist on keeping paired cusps, 
they also apply to the structure of the loops in the reduced states of this section. The basic reduction move in these figures corresponds to the 
elimination of two consecutive
cusps on a single loop by the move 3 of Figure 11 and also by move C of Figure 12. 
Note that we allow cancellation of consecutive cusps along a loop where the cusps both 
point to the same local region of the two local regions delineated by the loop. This is an abstraction of move 3 of Figure 11. We do not allow cancellation of a 
``zig-zag" where two consecutive cusps point to opposite (local) sides of the loop.  Note that in Figure 49, we have reformulated these loop reductions
in terms of special arrows attached tangentially to the loops. (Similar conventions are used in \cite{Miyazawa2}, and the reader should be warned that
our conventions are different from those in that paper.) Our convention is that {\it a cusp pointing upward and proceeding from left to right for the observer of the plane is labelled
with an arrow that goes from left to right.} We then see that each loop is labeled with arrows so that two consecutive arrows pointing in the same direction 
will cancel with each other. Opposed arrows do not cancel. Each loop gives rise (from the original state) to a sequence of arrows. The reduced loop
is coded by arrows where no further cancellation can be obtained. It is easy to see that the resulting cyclic sequence of arrows for the reduced loop
is unique and independent of the sequence of cancellations.  In the present section we use a complex reduction
procedure for the states, but the reduction of individual loops follows the same pattern as in Section 9. Here we pair the cusps  according
to the rules of Figures 11 and 12. Note that in a classical knot or link diagram, all state loops reduce to
loops that are free from cusps. See Proposition 2 in Section 9 for a proof of this fact.
\bigbreak 

\noindent {\bf Lemma 3.} {\it Each individual loop in a reduced state for the extended bracket is reduced according to the arrow rules of Section 9
 (as also discussed above).} 
\bigbreak

\noindent {\bf Proof.} 
Note that the type $C$ reduction, to which we have given precedence, involves an arrow reduction of the form 
$$[\cdots X \Leftarrow \Rightarrow \Rightarrow \Leftarrow Y \cdots] \sim [\cdots X \Leftarrow \Leftarrow Y \cdots] \sim [\cdots X Y \cdots]. $$
Here we refer to the arrow conventions of Figure 49, and we leave it to the reader to translate the arrow sequence for the type $C$ move.
In this formalism $ \Rightarrow \Rightarrow $ and $ \Leftarrow \Leftarrow $ correspond to two consecutive cusps that are not in zig-zag form. Such cusps cancel one
another in the reduction.
From the point of view of an indvidual loop, this reduction is just one choice among many reduction pathways for the unique reduction of the loop. For individual loops,
it is easy to see that the reduction is unique and independent of the order of reducing operations. To see this, regard $\alpha$ as $\Rightarrow$ and 
$ \beta $ as $\Leftarrow.$ Then we are using words in $\alpha$ and $\beta$ such that $\alpha^2$ and $\beta^2$ are empty words. Such words give a presentation
of the free group on $\alpha$ and $\beta$ modulo the relations that each has square equal to the identity. This is a presentation of the the free product of
$Z_{2}$ with $Z_{2},$ and the reduced words are alternating products of $\alpha$ and $\beta.$ For the state loops, the words are further reduced by cyclic
permutations so that the unique reduced forms are alternating products of $\alpha$ and $\beta$ that have an even number of terms (hence if they begin in $\alpha$ they
end in $\beta$ and if they begin in $\beta$ they end in $\alpha$).
Note that by our reduction rules, whenever there is a pair of consectutive cusps
along a loop with the same arrow assignment, there will be a reduction corresponding, at the loop level, to the cancellation of the arrows. Thus the loop ends up reduced in a canonical
fashion, while its cusp pairings are controlled by the rules given in this section. This completes the proof of the Lemma. //
\bigbreak

\begin{figure}
     \begin{center}
     \begin{tabular}{c}
     \includegraphics[width=10cm]{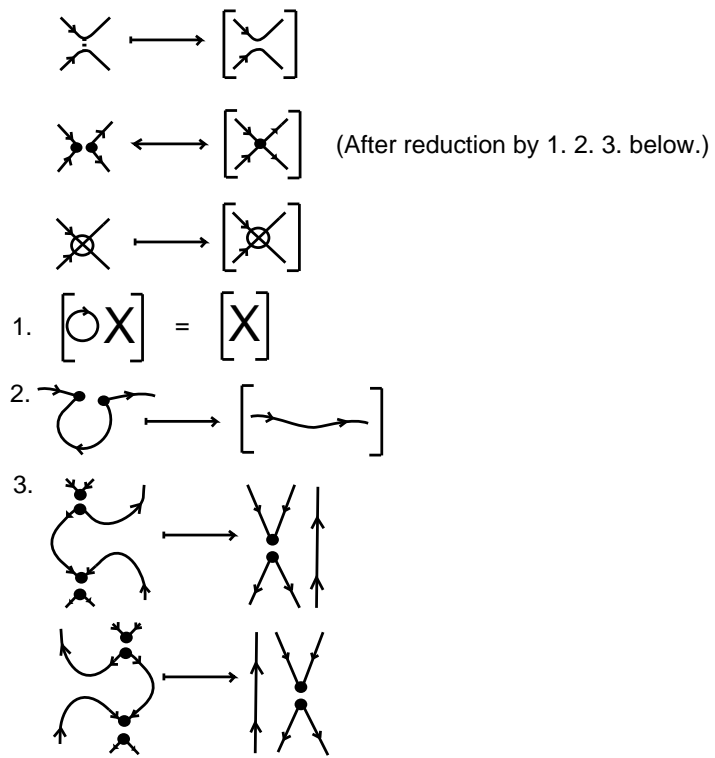}
     \end{tabular}
     \caption{\bf Basic Replacements}
     \label{Figure 11}
\end{center}
\end{figure}

\begin{figure}
     \begin{center}
     \begin{tabular}{c}
     \includegraphics[width=4cm]{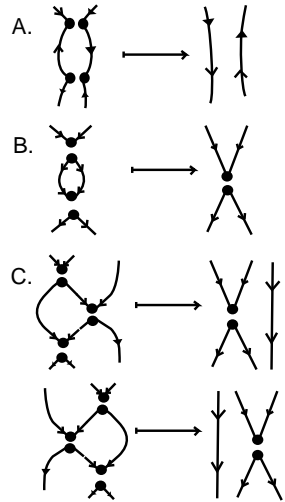}
     \end{tabular}
     \caption{\bf Special Replacements}
     \label{Figure 12}
\end{center}
\end{figure}

\begin{figure}
     \begin{center}
     \begin{tabular}{c}
     \includegraphics[width=7cm]{F13.EPSF}
     \end{tabular}
     \caption{\bf The Type One Move}
     \label{Figure 13}
\end{center}
\end{figure}

\begin{figure}
     \begin{center}
     \begin{tabular}{c}
     \includegraphics[width=7cm]{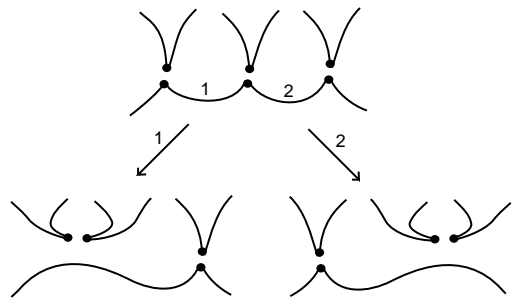}
     \end{tabular}
     \caption{\bf Multiplicity}
     \label{Figure 14}
\end{center}
\end{figure}

\begin{figure}
     \begin{center}
     \begin{tabular}{c}
     \includegraphics[width=9cm]{F15.EPSF}
     \end{tabular}
     \caption{\bf Uniqueness of Special Replacement $A$}
     \label{Figure 15}
\end{center}
\end{figure}

\begin{figure}
     \begin{center}
     \begin{tabular}{c}
     \includegraphics[width=6cm]{F16.EPSF}
     \end{tabular}
     \caption{\bf Uniqueness of Special Replacement $A$}
     \label{Figure 16}
\end{center}
\end{figure}

\begin{figure}
     \begin{center}
     \begin{tabular}{c}
     \includegraphics[width=6cm]{F17.EPSF}
     \end{tabular}
     \caption{\bf Uniqueness of Special Replacement $A$}
     \label{Figure 17}
\end{center}
\end{figure}

\begin{figure}
     \begin{center}
     \begin{tabular}{c}
     \includegraphics[width=4cm]{F18.EPSF}
     \end{tabular}
     \caption{\bf Uniqueness of Special Replacement $B$}
     \label{Figure 18}
\end{center}
\end{figure}

\begin{figure}
     \begin{center}
     \begin{tabular}{c}
     \includegraphics[width=6cm]{F19.EPSF}
     \end{tabular}
     \caption{\bf Uniqueness of Special Replacement $B$}
     \label{Figure 19}
\end{center}
\end{figure}

\begin{figure}
     \begin{center}
     \begin{tabular}{c}
     \includegraphics[width=6cm]{F20.EPSF}
     \end{tabular}
     \caption{\bf Uniqueness of Special Replacement $B$}
     \label{Figure 20}
\end{center}
\end{figure}

\begin{figure}
     \begin{center}
     \begin{tabular}{c}
     \includegraphics[width=7cm]{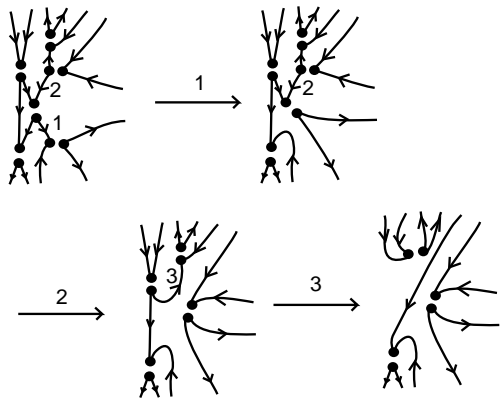}
     \end{tabular}
     \caption{\bf Special Replacement $C$ Requires a Precedence Rule}
     \label{Figure 21}
\end{center}
\end{figure}

\begin{figure}
     \begin{center}
     \begin{tabular}{c}
     \includegraphics[width=6cm]{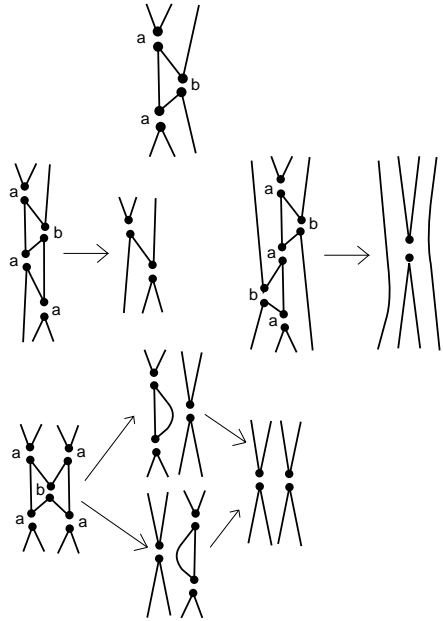}
     \end{tabular}
     \caption{\bf Networks of $C$ - Moves}
     \label{Figure 22}
\end{center}
\end{figure}

\begin{figure}
     \begin{center}
     \begin{tabular}{c}
     \includegraphics[width=9cm]{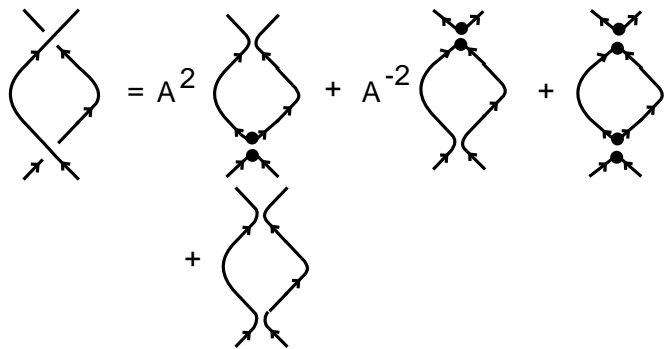}
     \end{tabular}
     \caption{\bf Oriented Second Reidemeister Move}
     \label{Figure 23}
\end{center}
\end{figure}

\begin{figure}
     \begin{center}
     \begin{tabular}{c}
     \includegraphics[width=8cm]{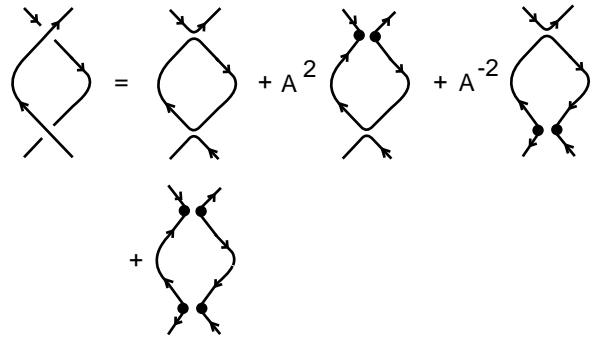}
     \end{tabular}
     \caption{\bf Reverse Oriented Second Reidemeister Move}
     \label{Figure 24}
\end{center}
\end{figure}

\begin{figure}
     \begin{center}
     \begin{tabular}{c}
     \includegraphics[width=10cm]{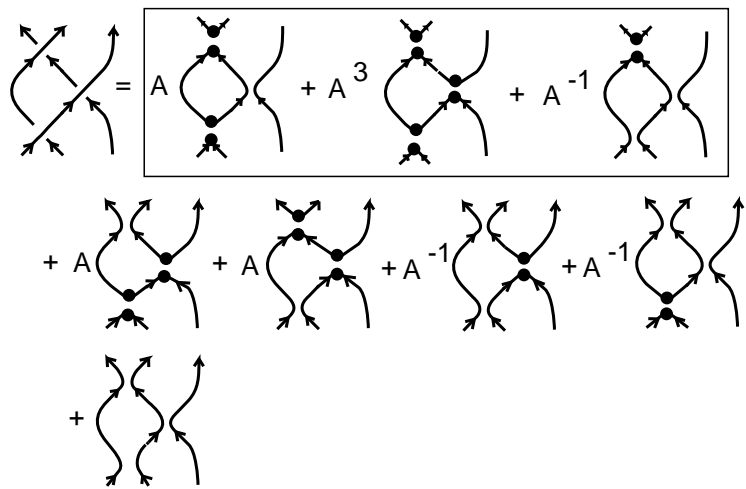}
     \end{tabular}
     \caption{\bf Third Reidemeister Move}
     \label{Figure 25}
\end{center}
\end{figure}

\begin{figure}
     \begin{center}
     \begin{tabular}{c}
     \includegraphics[width=8cm]{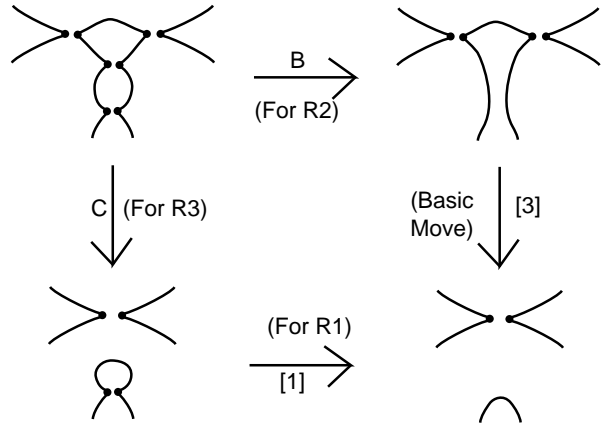}
     \end{tabular}
     \caption{\bf Reduction Relation}
     \label{Figure 26}
\end{center}
\end{figure}

Here is a first example of a calculation of the extended bracket invariant. View Figure 27. The virtual knot $K$ in this figure has two crossings.
One can see that this knot is a non-trival virtual knot by simply calculating the odd writhe $J(K)$ (defined in section 5). We have that
$J(K) = 2,$ proving that $K$ is non-trivial and non-classical. This is the simplest virtual knot, the analog of the trefoil knot for virtual knot 
theory. The extended bracket polynomial gives an independent verification that $K$ is non-trivial and non-classical.
Note that in the calculation of $<<K>>$ the fourth state entails a special replacement, and in the final polynomial we have a monomial term and
a term involving a non-trivial flat class of the simple flat $H$ (shown in the figure) consisting of two circles with one real flat crossing and one virtual
crossing. 
\bigbreak

In this replacement, the reader should note that there are two loops, one without self-crossings and one with a virtual self-crossing. The loop
with the virtual self-crossing does not fall under the category of our reduction rules. It is the loop with no self-crossing to which the reduction rules 
apply. The reduction rules are sensitve to the order of edges at each disoriented site. 
\bigbreak

In Figure 27 the virtual graph $H$ is non-trivial because any flat link diagram of two components  with an odd number of real crossings between the components is
necessarily non-trivial (since this parity is preserved by the Reidemeister moves). The appearance of the non-trivial virtual graph $H$ shows that $K$ is
non-classical. 
\bigbreak

\begin{figure}
     \begin{center}
     \begin{tabular}{c}
     \includegraphics[width=10cm]{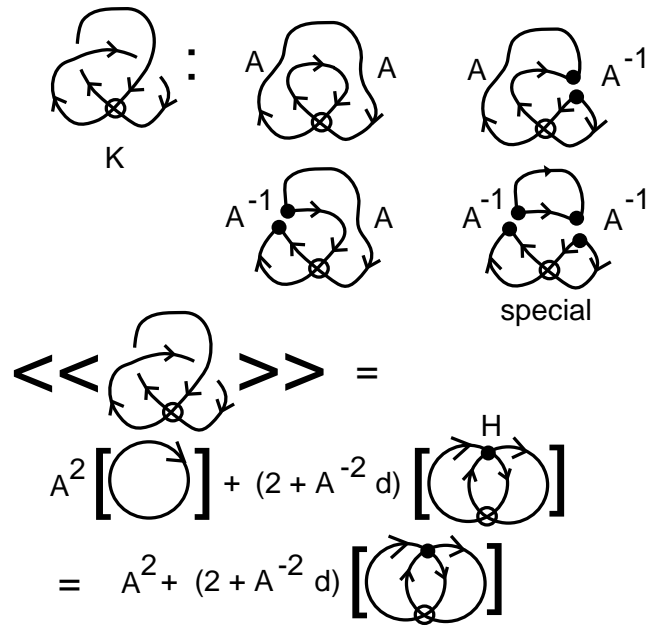}
     \end{tabular}
     \caption{\bf Example1}
     \label{Figure 27}
\end{center}
\end{figure}

In the next example, shown in Figure 28, we have a long virtual diagram $L$ with two crossings. The calculation of the extended bracket for 
$L$ is given in Figure 28 and shows that it is a non-trivial and non-classical. In fact, this same formalism
proves that $Flat(L)$ is a non-trivial flat link (as we have also seen from Section 5.). {\it Note that the extended bracket is an invariant of 
flat diagrams when we take $A=1.$} With $A=1$ we have $d=-2$ and so $<<FL>> = [\lambda]+2[OH]+[\lambda G]$ where the graph labels in this formula
are shown in Figure 28. We will later give other examples of the detection of flat
non-triviality that are more complex than this example.
\bigbreak

Following Example 2, we have the calculation in Figure 29 of the closure of $L$ forming the virtual diagram $CL.$ This diagram is unknotted, and 
the calculation reflects this fact. Note that in this calculation there appears a special replacement, corresponding to a Reidemeister 2 move, that is not
available in the correspondiing long diagram. Thus we see that the extended bracket can (sometimes) discriminate between a long knot and its closure.
\bigbreak

\begin{figure}
     \begin{center}
     \begin{tabular}{c}
     \includegraphics[width=10cm]{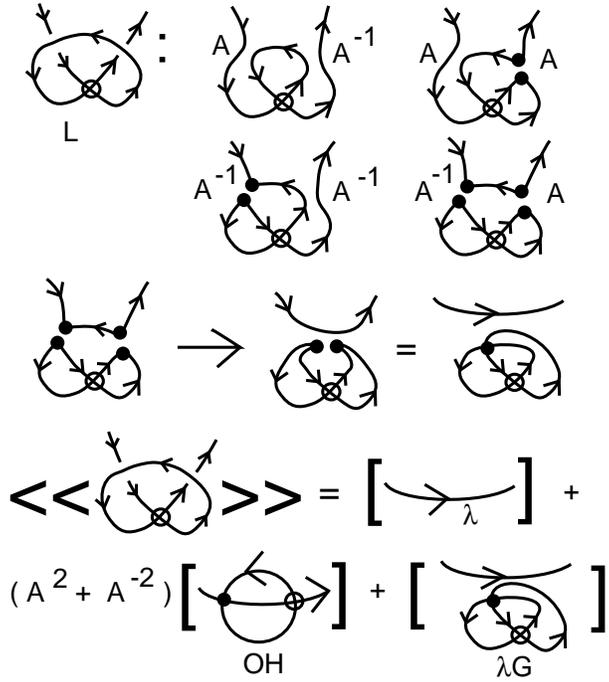}
     \end{tabular}
     \caption{\bf Example2}
     \label{Figure 28}
\end{center}
\end{figure}

\begin{figure}
     \begin{center}
     \begin{tabular}{c}
     \includegraphics[width=10cm]{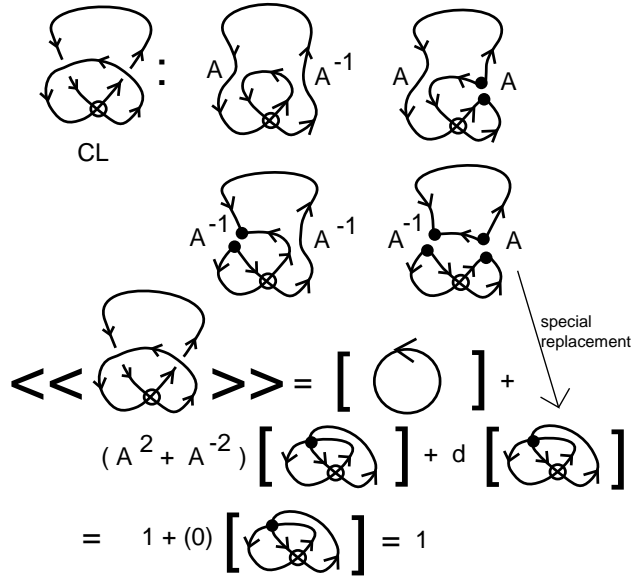}
     \end{tabular}
     \caption{\bf Example2.1}
     \label{Figure 29}
\end{center}
\end{figure}

\noindent The application of the special replacements requires some attention to context. In Figures 30 and 31 we
illustrate cases where there are states with special replacement disoriented loops that are complicated by the presence of virtual  crossings
(eliminated by a detour move) and local curl reversals. At the Gauss diagram level the virtual crossings are simply not present, and they must be discounted 
in searching for special replacements. In Figure 31 we give an example of a calculation of the extended bracket on a long flat virtual that is in fact
trivial (we leave the verification of its triviality to the reader). Note that in the course of the calculation we find a single disoriented loop that can be
removed and this loop has virtual intersections with the rest of the diagram while the vertex orders of its site interactions correspond to the basic
reduction rules. It is removed and the resulting flat diagrams collect and  cancel. 
\bigbreak

The next example is given in Figures 32, 33 and 34. Here we calculate the extended bracket for a non-trivial virtual knot with unit Jones polynomial.
The reader can follow the calculation through these figures and see that the extended bracket detects the non-triviality of this knot and proves that it
is not isotopic to a classical knot diagram. In this example we find that the extended bracket is the sum of a scalar term and a polyomial multiple of 
a virtual graph  $F.$  A proof that the flat virtual corresponding to the virtual graph $F$ is non-trivial is shown
in Figure 35. We label arcs in the diagram and prolong these labels through the virtual crossings, but multiply them by a commuting variable $s$ or $s^{-1}$
depending upon a right of left passage as shown in this figure. The resulting  module structure over the ring $Z[s,s^{-1}]$ is an invariant of the flat
virtual diagram. This module is an example of a biquandle invariant for flat virtuals. More information and definitions for biquandles can be found in
\cite{KM1,KM2,FJK}. For our purposes, we only need to know that the virtual graph $F$ is non-trivial, and this can be seen directly via the cyclic ordering of 
the edges at the vertices and an application of the Jordan curve theorem to the cycles in the graph. We leave the details of this approach to the reader.
\bigbreak

Figures 36, 37, 38, 39, 40 and 41 exhibit the calculation of the extended bracket for the Kishino diagram. Figure 36 is a list of all the states of the Kishino
diagram. Figure 37 is drawn in parallel to Figure 36 and makes reductions on those states that can be reduced uniquely. Boxes in Figure 37 are placed around those states
that have a multiple reduction. In going from Figure 37 to Figure 38, further reductions are made. For example, the state in Figure 37 just above the lower left-hand
corner reduces to a simple loop and so can be eliminated. In Figure 38 this position appears with just the coefficient $A^{-4}.$ All states in Figure 38 that are
not in boxes are fully reduced. Figures 39 and 40 illustrate the two basic cases that involve multiplicity in the reduction. Finally Figure 41 gives the 
final result of the calculation.  Note that in the next to final result of Figure 41 we have left the states in reduced form, but have not
replaced the paired cusps with graphical nodes.
As we have remarked earlier, the replacement by graphical nodes is really a notational device to indicate the end of the reduction process. In this example
we leave the result in the form of reduced states so that the reader can compare it with the arrow polynomial (simple extended bracket) to be introduced
in Section 9. In that formulation, each reduced state is replaced by a certain monomial. We give the graphical form of the extended bracket for the Kishino knot
in the last equality of Figure 41.
\bigbreak

In this example we prove that the Kishino diagram is
non-trivial, and in fact, by taking $A=1$, we see that the flat Kishino diagram is non-trivial. In the case of the Kishino diagram, this adds one more proof
to a long list of verifications of its non-triviality. See the problem list of Fenn, the author and Manturov \cite{VP}. It is of interest to see how the
extended bracket invariant sees the non-triviality of the Kishino flat diagram as shown in the summary result of Figure 41.
\bigbreak

We verify that the flat diagram shown in Figure 42 is non-trivial by using the extended bracket. The results of the state reductions are shown in the 
small diagrams in that figure.  
\bigbreak

The last example is shown in Figures 43 and Figure 44. In these figures we show the result of expanding a virtualized
classical crossing using the extended bracket state sum. Virtualization of a crossing was described in Section 6. In a
virtualized crossing, one sees a classical crossing that is flanked by two virtual crossings. In Section 6 we showed that
the standard bracket state sum does not see virtualization in the sense that it  has the same value as the result of
smoothing both flanking virtual crossings that have been added to the diagram. The result is that the value of of the
bracket polynomial of the knot with a virtualized classical crossing is the same as the value of the bracket polynomial of
the original knot after the same crossing has been {\it switched} (exchanging over and undercrossing segments). 
\bigbreak

As one can see from the formula in Figure 43, this smoothing property of the bracket polynomial will not generally be the case for the extended bracket state
sum. In Figure 44 we show that this difference is indeed the case for an infinite collection of examples. In that figure we use a tangle $T$ that is assumed
to be a classical tangle. Extended bracket expansison of this tangle is necessarily of the form shown in that figure: a linear combination of an oriented
smoothing and a reverse oriented smoothing with respective coefficients $a(T)$ and $b(T)$ in the Laurent polynomial ring $Q[A,A^{-1}].$ We leave the
verification of this fact to the reader. In Figure 44 we show a generic diagram that is obtained by a {\it single} virtualization from a classical diagram,
and we illustrate the calculation of its extended bracket invariant. As the reader can see from this Figure, there is a non-trivial graphical term whenever
$b(T)$ is non-zero. Thus we conclude that the single virtualization of any classical link diagram (in the form shown in this Figure) will be non-trivial and
non-classical whenever $b(T)$ in non-zero. This is an infinite class of examples, and the result can be used to recover the results about single
virtualization that we  obtained in a previous paper with Heather Dye \cite{MinSurf} using the surface bracket polynomial.
\bigbreak

\begin{figure}
     \begin{center}
     \begin{tabular}{c}
     \includegraphics[width=6cm]{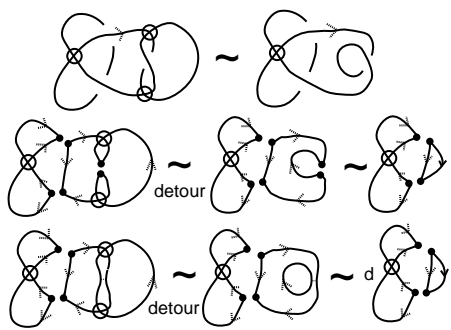}
     \end{tabular}
     \caption{\bf Effect of the Detour Move}
     \label{Figure 30}
\end{center}
\end{figure}

\begin{figure}
     \begin{center}
     \begin{tabular}{c}
     \includegraphics[width=8cm]{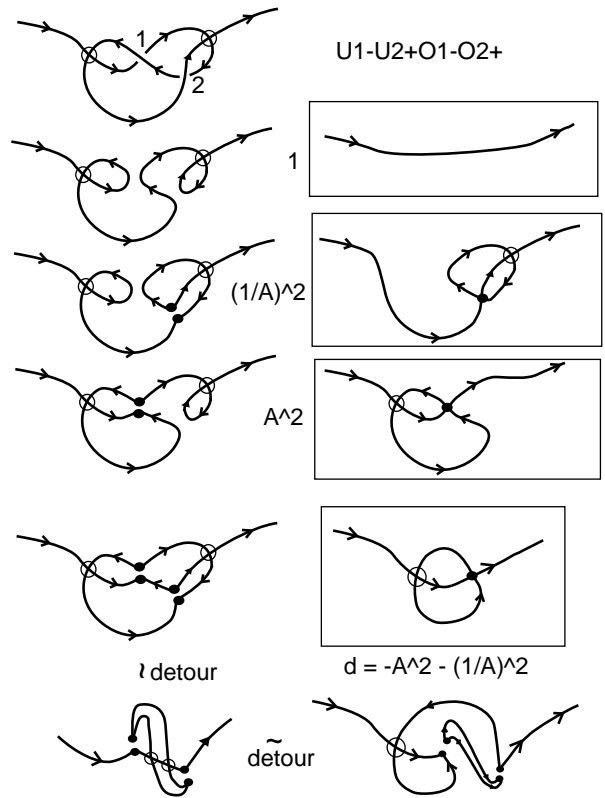}
     \end{tabular}
     \caption{\bf Cancellation in a Trivial Long Virtual}
     \label{Figure 31}
\end{center}
\end{figure}

\begin{figure}
     \begin{center}
     \begin{tabular}{c}
     \includegraphics[width=8cm]{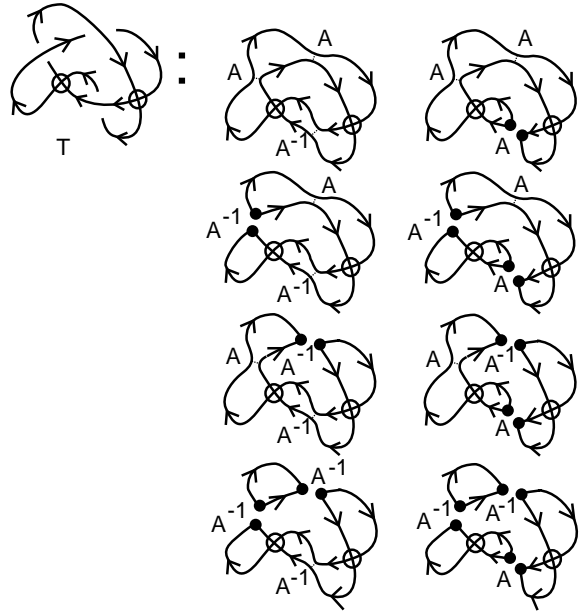}
     \end{tabular}
     \caption{\bf Virtualized Trefoil States}
     \label{Figure 32}
\end{center}
\end{figure}

\begin{figure}
     \begin{center}
     \begin{tabular}{c}
     \includegraphics[width=10cm]{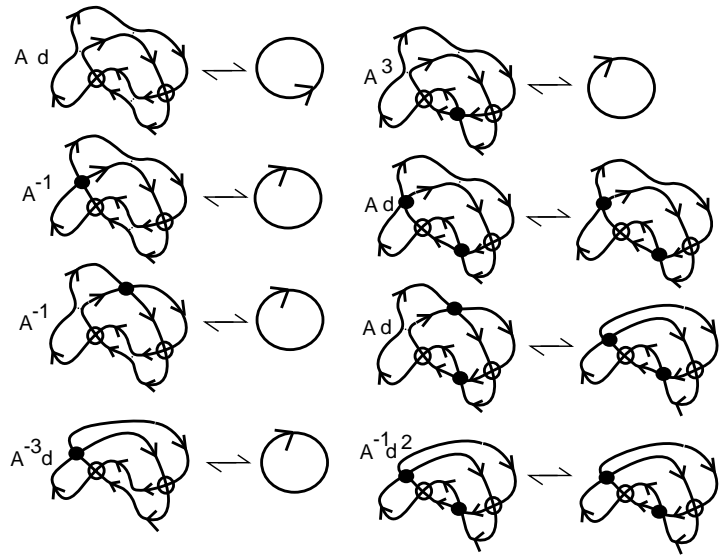}
     \end{tabular}
     \caption{\bf Flattened Virtualized Trefoil States}
     \label{Figure 33}
\end{center}
\end{figure}

\begin{figure}
     \begin{center}
     \begin{tabular}{c}
     \includegraphics[width=7cm]{F34.EPSF}
     \end{tabular}
     \caption{\bf Extended Bracket for the Virtualized Trefoil}
     \label{Figure 34}
\end{center}
\end{figure}

\begin{figure}
     \begin{center}
     \begin{tabular}{c}
     \includegraphics[width=5cm]{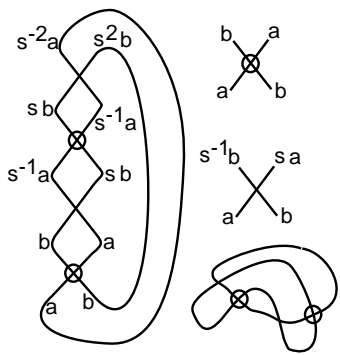}
     \end{tabular}
     \caption{\bf Verification of Non-triviality of a Specific Flat Link}
     \label{Figure 35}
\end{center}
\end{figure}

\begin{figure}
     \begin{center}
     \begin{tabular}{c}
     \includegraphics[width=7cm]{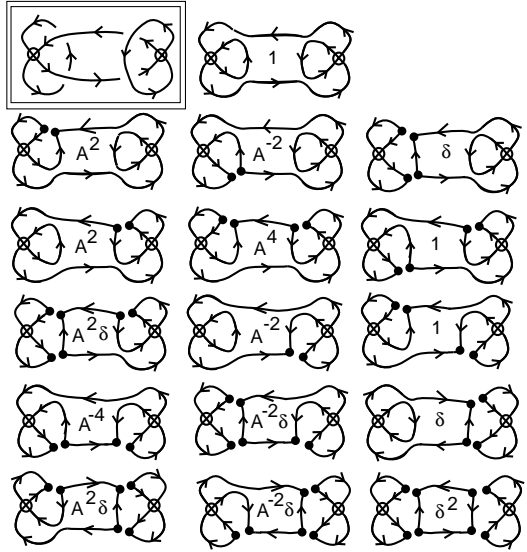}
     \end{tabular}
     \caption{\bf Kishino Diagram States}
     \label{Figure 36}
\end{center}
\end{figure}

\begin{figure}
     \begin{center}
     \begin{tabular}{c}
     \includegraphics[width=8cm]{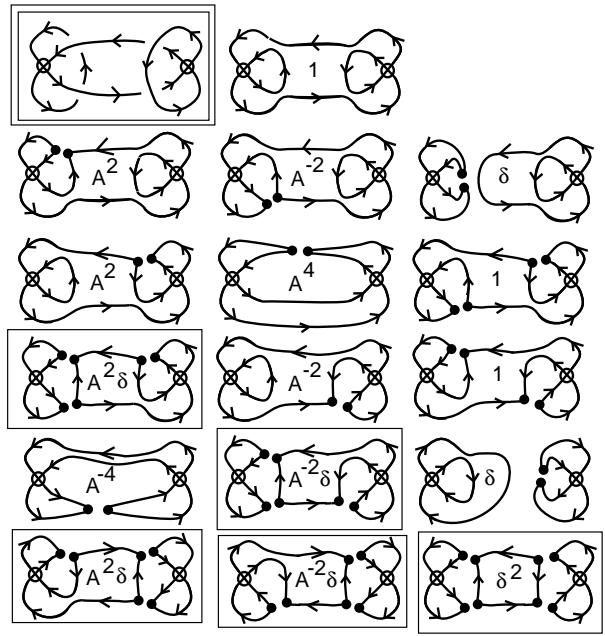}
     \end{tabular}
     \caption{\bf Reducing Kishino Diagram States}
     \label{Figure 37}
\end{center}
\end{figure}

\clearpage

\begin{figure}
     \begin{center}
     \begin{tabular}{c}
     \includegraphics[width=7cm]{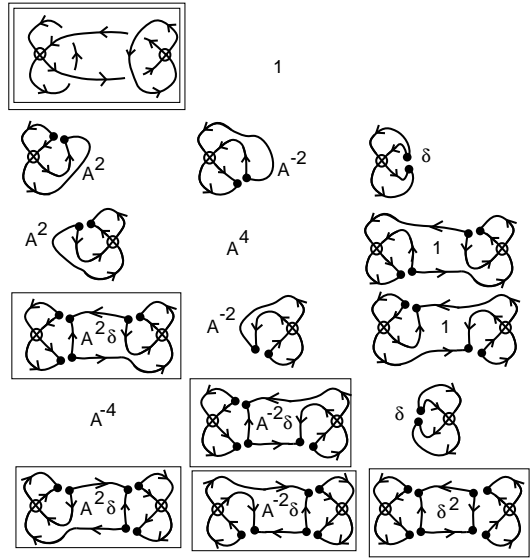}
     \end{tabular}
     \caption{\bf Reducing Kishino Diagram States}
     \label{Figure 38}
\end{center}
\end{figure}

\begin{figure}
     \begin{center}
     \begin{tabular}{c}
     \includegraphics[width=6cm]{F39.EPSF}
     \end{tabular}
     \caption{\bf Reducing in Multiplicity}
     \label{Figure 39}
\end{center}
\end{figure}

\begin{figure}
     \begin{center}
     \begin{tabular}{c}
     \includegraphics[width=8cm]{F40.EPSF}
     \end{tabular}
     \caption{\bf Reducing in Multiplicity}
     \label{Figure 40}
\end{center}
\end{figure}

\begin{figure}
     \begin{center}
     \begin{tabular}{c}
     \includegraphics[width=8cm]{F41.EPSF}
     \end{tabular}
     \caption{\bf Extended Bracket for the Kishino Diagram}
     \label{Figure 41}
\end{center}
\end{figure}

\begin{figure}
     \begin{center}
     \begin{tabular}{c}
     \includegraphics[width=8cm]{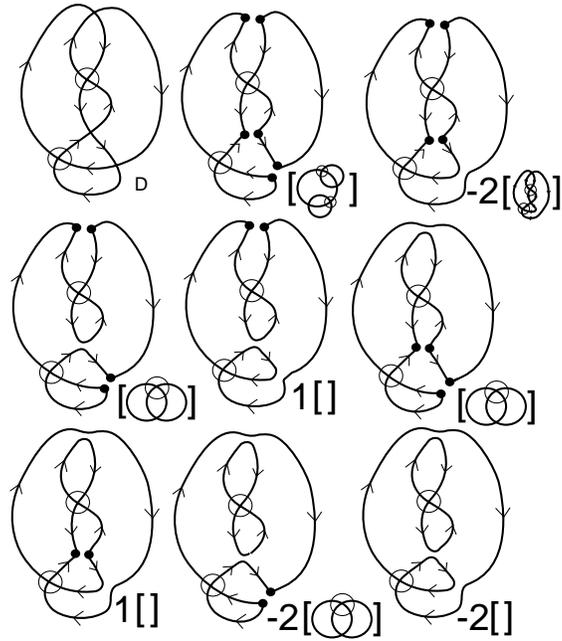}
     \end{tabular}
     \caption{\bf Verifying the Non-Triviality of a Reclacitrant Flat Diagram}
     \label{Figure 42}
\end{center}
\end{figure}

\begin{figure}
     \begin{center}
     \begin{tabular}{c}
     \includegraphics[width=8cm]{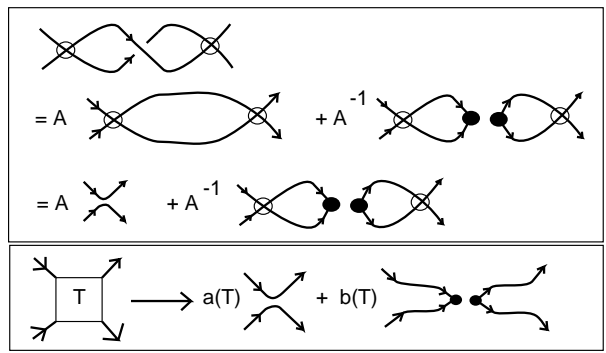}
     \end{tabular}
     \caption{\bf The $1$-Virtualization of a Classical Diagram (A)}
     \label{Figure 43}
\end{center}
\end{figure}

\begin{figure}
     \begin{center}
     \begin{tabular}{c}
     \includegraphics[width=10cm]{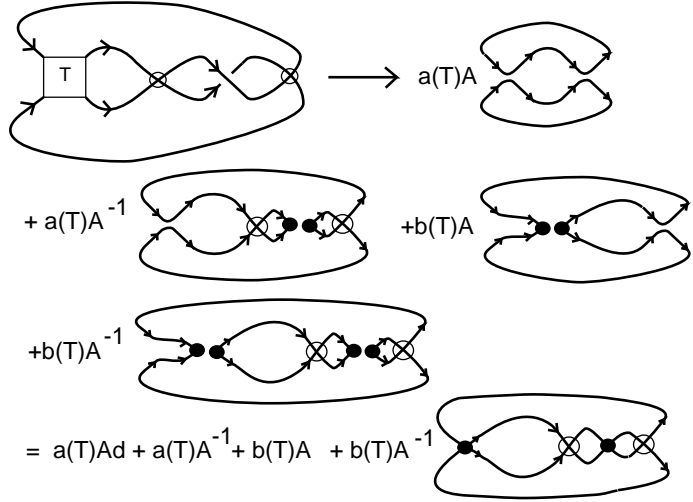}
     \end{tabular}
     \caption{\bf The $1$-Virtualization of a Classical Diagram (B)}
     \label{Figure 44}
\end{center}
\end{figure}
\bigbreak
 
\section{Estimating Virtual Crossing Number}
The {\it virtual crossing number} $VC(K)$ of a virtual knot or
link $K$ is the least number of virtual crossings in any diagram
that represents the isotopy class of $K.$ In computing the 
extended bracket polynnomial, we create a collection of reduced
virtual graphs $G$ from the states of a diagram $K$. The
virtual isotopy class of each graph receiving  a non-zero
coefficient in $<<K>>$ is an isotopy invariant of the orginal
virtual knot or link $K.$ Each graph also has a virtual crossing
number, that we shall denote by $vc(G).$ We then have the
following theorem.
\bigbreak

\noindent {\bf Theorem 5.} Let $K$ be a virtual knot or link diagram
and let $vc(K)$ denote the maximum value of $vc(G)$ over the virtual
crossing numbers for the graphs $G$ that appear with non-zero
coefficients in the extended bracket state sum $<<K>>.$
Then the virtual crossing number for $K$ is bounded below by
$vc(K).$ That is, 
$VC(K) \ge vc(K).$ 
\bigbreak

\noindent {\bf Proof.} Since the graphs with non-zero coefficients in $<<K>>$ are individually invariants
of $K$, it follows that their individual virtual crossing numbers are invariants of $K$. Furthermore,
each such individual crossing number is necessarily less than the virtual crossing number of $K$, since 
the extended bracket of $K$ can be computed from a representative of $K$ with minimal virtual crossing number
(and the resulting graphs will all have no more than this minimal number of virtual crossings). This completes
the proof of the Theorem. //
\bigbreak

We can use this Theorem directly to estimate, and sometimes compute, the virtual crossing number of specific examples.
Thus we conclude that the virtual crossing number for Example 1 of Figure 27 is one. The virtual crossing number for
the long knot of Example 2 of Figure 28 is one. The virtual crossing number of the virtualized trefoil knot of Figure
33 is two. Note that in this case, the genus (least genus surface on which the knot is represented) is one. For the Kishino
diagram, the calculation summed up in Figure 41 shows that it has virtual crossing number two (and it also has genus two
via our work on the surface bracket polynomial \cite{MinSurf}). The flat virtual of Figure 42 has virtual crossing number 
equal to two. 
\bigbreak

We can combine this estimation theorem with other facts about the extended bracket polynomial to get estimates without 
actually computing the entire state sum. In order to accomplish this, we generalize a well-known result about the 
bracket polynomial to virtual knots and links. For this purpose, we call a site in a state $S$ of $K$ {\it self-touching}
if both local edges of this site belong to the same loop in $S.$
\bigbreak

\begin{figure}
     \begin{center}
     \begin{tabular}{c}
     \includegraphics[width=6cm]{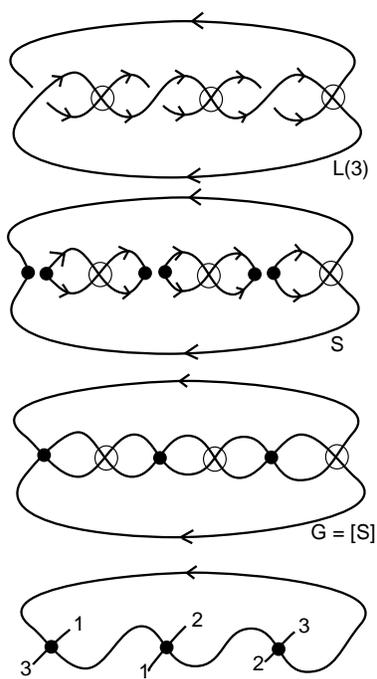}
     \end{tabular}
     \caption{\bf $L(n)$ Has Genus One and Virtual Crossing Number $n$.}
     \label{Figure 45}
\end{center}
\end{figure}

\begin{figure}
     \begin{center}
     \begin{tabular}{c}
     \includegraphics[width=6cm]{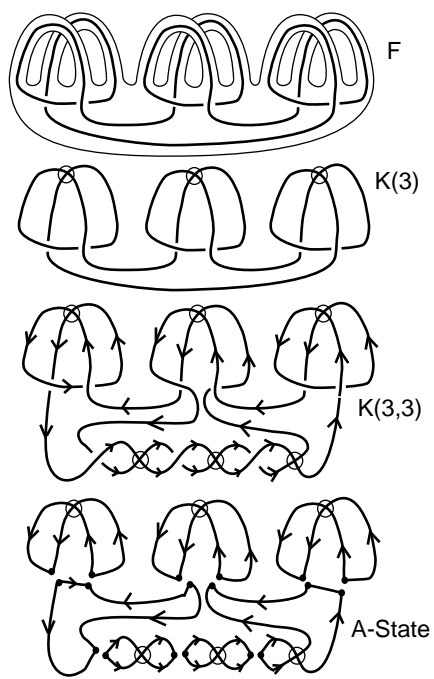}
     \end{tabular}
     \caption{\bf $K(n,m)$ has genus $n+1$ and Virtual Crossing Number $n+m$.}
     \label{Figure 46}
\end{center}
\end{figure}

\noindent {\bf Proposition 1.} Let $K$ be a virtual knot or link diagram. Let $S$ be the $A$-state of $K$ and $S'$ be the 
$A^{-1}$-state of $K$ (all smoothings of type $A$ or all smoothings of type $A^{-1}$ respectively). If no
site in $S$ is self-touching, then the maximal degree of $<K>$ is attained on $S$ and it is $v(K) + 2||S|| - 2$
where $v(K)$ is the number of crossings in the diagram $K$ and $||S||$ is the number of loops in the state $S.$
Similarly, if no site in $S'$ is self-touching, then the minimal degree of $<K>$ is attained on $S'$ and it is 
$-v(K) - 2||S'|| +2.$ Since the extended bracket state sum $<<K>>$ is a refinement of the bracket state sum, we conclude
that the reduced virtual graphs $[S]$ and $[S']$ will have non-zero coefficients in the 
extended bracket. Hence, if $S$ is not self-touching or if $S'$ is not self-touching, then $VC(K) \ge vc([S])$ or $VC(K) \ge vc([S']).$
\bigbreak

\noindent {\bf Proof.} For the degree inequalities stated in the first part of the Proposition, see the author's 
papers and books \cite{K,KNOTS} and the book by Lickorish \cite{Lickorish}. Knots and links that satisfy the
non-self-touching conditions for both $S$ and $S'$ are said to be {\it adequaate}. Knots and links that satisfy one but not both of these conditions are said to be
{\it semi-adequate}. The rest of the statement follows at once from our work on the extended bracket state sum and the previous Theorem. //
\bigbreak

For an example of the application of this proposition, view Figure 45. In that figure we illustrate a virtual link
$L(3)$, part of a class of virtual links $L(n)$ for $n= 1,2,3, \cdots.$ The $A$-state $S$ of $L(n)$ reduces to a graph
with virtual crossing number $N$ and the $A$-state has no self-touching sites. It follows that the virtual crossing number
of $L(n)$ is equal to $n.$ It is easy to see that the genus of $L(n)$ is one. The figure for $L(3)$ illustrates the
$A$-state $S$ and the resulting reduced graph $G = [S].$ The figure also illustrates how the graph consists of two
intersecting curves, each a Jordan curve in the plane. The connections that one curve must make in relation to the other
requires three virtual crossings for $L(3).$ Similarly, $n$ virtual crossings are required for $L(n).$
\bigbreak

The next example is shown in Figure 46. This figure illustrates a knot $K(n)$, for $n=3$, by first showing it as a diagram on
a surface $F$ with boundary of genus $g = n = 3$, and then showing the virtual diagram for $K(n)$ as a projection from this abstract link diagram.
One can apply the surface bracket \cite{MinSurf} to the knot in $F$ to show that $K(n)$ has minimal genus $n.$ We omit the details of this calculation, but
point out that the reader can verify it by showing that the states of the surface bracket expansion of $K(n)$ span the first homology group of the surface $F.$
Below the diagram for $K(n)$ we have indicated the diagram for $K(n,m)$, obtained from $K(n)$ by adding a two-strand virtual braid with $m$ real crossings and 
$m$ virtual crossings. In the figure $n=m=3.$ A similar calculation with the surface bracket shows that $K(n,m)$ has minimal genus $n+1$ (the two-braid is
supported on an extra handle) for $m > 0.$ Returning to our extended bracket invariant, we note that the $A$-state $S$ for $K(n,m)$ has no self-touching loops,
and its reduced graph has $n+m$ virtual crossings (we leave this for the reader to check). Thus $K(n,m)$ has minimal genus $n+1$ and virtual crossing number
of its reduction is $n+m.$ Thus, by a combination of the surface bracket and the extended bracket state sum, we have shown that there exist virtual knots $K$ with arbitrarily 
high minimal genus and arbitrarily high virtual crossing number with the gap $VC(K)-g(K)$ positive and arbitrarily large.
\bigbreak   

Refer to Figure $47.$ Here we give an example of a virtual knot that is undetectable by the extended bracket.
The reader can check that this knot is also undectected by the Sawollek polynonmial \cite{Sawollek} but that it can be seen
to be non-trivial and non-classical by examining the structure of its Alexander module (Its classical Alexander polynomial is non-trivial and 
does not have the symmetry property of a classical knot). This example (for which we are indebted to Slavik Jablan) shows that the extended bracket 
invariant is not invulnerable. In joint work with Slavik Jablan \cite{Slavik}, we will publish tables of calculations for the simple extended bracket
and its relatives. To see why this example is not dectectable by the extended bracket polynomial, examine the rest of Figure 48, which shows that 
the polynomial is invariant under the  ``double replacement" illustrated there. This replacement is similar to the virtualization replacement
discussed earlier in this paper. It is equivalent to the ``double flype" studied in \cite{Kamada1,Kamada2}. It is easy to see that 
the example in Figure 47 is unknotted after a double replacement. In the figure, the undetectable knot is labeled $KS$ and we show that it contains a
configuration lableled $S$ in that figure. The corresponding configuration $T$ can replace $S$ without changing the extended bracket, as is shown in the
figure. We illustrate  the replacement for $KS$ with an arrow to the new knot $KT$. It is easy to see that $KT$ is unknotted.
\bigbreak

\begin{figure}
     \begin{center}
     \begin{tabular}{c}
     \includegraphics[width=8cm]{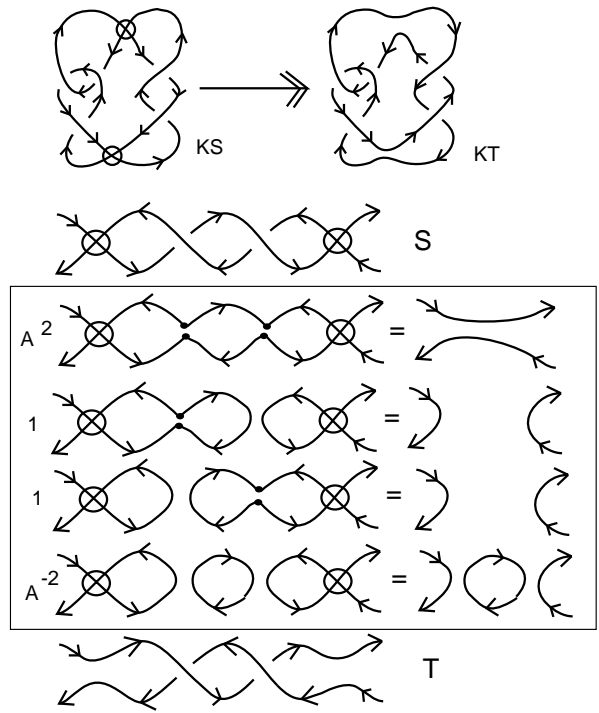}
     \end{tabular}
     \caption{\bf A Virtual Knot Undetectable by the Extended Bracket.}
     \label{Figure 47}
\end{center}
\end{figure}

\section{The Arrow Polynomial}
In this section we describe and compute a simplification of the extended bracket invariant described in this paper.
This simple extended bracket invariant will be called the {\it arrow polynomial}. In \cite{DyeKauff} we investigate the arrow polynomial more deeply
and show how it can be used to estimate virtual crossing numbers.
\bigbreak

In the previous sections we have reduced the states until there is a remaining set of disoriented smoothings. We then welded the 
remaining disoriented smoothings to form graphs. The extended invariant is a sum of graphs (taken up to virtual equivalence in the plane)
weighted by polynomials. If, instead of welding these sites, we release the remaining smoothings we will get simpler graphs composed of disjoint collections of 
circle graphs that are labelled with the orientation markers and left-right distinctions that occur in the state expansion.  The basic conventions for this 
simplification are shown in Figures 48 and 49. In those figures we illustrate how the disoriented smoothing is a local disjoint union of two vertices.
Each vertex is denoted by an angle with arrows either both entering the vertex or both leaving the vertex. Furthermore, the angle locally divides the 
plane into two parts: One part is the span of an acute angle (of size less than $\pi$); the other part is the span of an obtuse angle. We refer to the span of the acute angle as
the {\it inside} of the vertex. In Figure 48 we have labeled the insides of the vertices with the symbol $\sharp .$ In making drawings of these graphs, one can
maintain the distinction between inside and outside of the vertices without the necessity of drawing the $\sharp.$ The reader will note that we have
consistently maintained the convention in the earlier parts of the paper where the vertices of a smoothing were always paired. Now it is essential to 
articulate this convention since we shall allow the pairs to come apart.
\bigbreak

Figure 48 illustates the basic reduction rule for the {\it arrow polynomial} ( or {\it simple extended bracket polynomial}) We shall denote this invariant by the notation
${\cal A}[K],$ for a virtual knot or link diagram $K.$ The reduction rule allows the cancellation of two adjacent vertices when they have {\it insides on the same
side} of the segment that connects them. When the insides of the vertices are on opposite sides of the connecting segment, then no cancellation is allowed.
All graphs are taken up to virtual equivalence, as explained earlier in this paper. Figure 48 illustrates the simplification of two circle graphs. In one case
the graph reduces to a circle with no vertices. In the other case there is no further cancellation, but the graph is equivalent to one without a virtual crossing.
The state expansion for ${\cal A}[K]$ is exactly as shown in Figure 10, but we use the reduction rule of Figure 48 so that each state is a disjoint union of reduced circle 
graphs. Since such graphs are planar, each is equivalent to an embedded graph (no virtual crossings) and the reduced forms of such graphs have $2n$ vertices that 
alternate in type around the circle so that $n$ are pointing inward and $n$ are pointing outward. The circle with no vertices is evaluated as $d = -A^2 - A^{-2}$ as
is usual for these expansions and the circle is removed from the graphical expansion. Let $K_{n}$ denote the circle graph with $2n$ alternating vertex types as shown in
Figure 48 for $n=1$ and $n=2.$
By our conventions for the extended bracket polynomial, each circle graph contributes $d = -A^2 - A^{-2}$ to the state sum and the graphs $K_{n}$ for $n \ge 1$ remain
in the graphical expansion. For the simplified version ${\cal A}[K]$ we can regard each $K_{n}$ as an extra variable in the polynomial. Thus a product of the $K_{n}$'s
denotes a state that is a disjoint union of copies of these circle graphs with multiplicities. By evaluating each circle graph as $d = -A^2 - A^{-2}$ we guarantee
that the resulting polynomial will reduce to the original bracket polynomial when each of the new variables $K_{n}$ is set equal to unity. Note that we 
continue to use the caveat that an isolated circle or circle graph (i.e. a state consisting in a single circle or single circle graph) is assigned a loop value
of unity in the state sum. This assures that ${\cal A}[K]$ is normalized so that the unknot receives the value one.
\bigbreak

In Figure 49 we show the {\it arrow convention} for rewriting a state loop in cusp form as a loop with a seqence of arrows, one arrow for each cusp.
The convention for replacing a tangential arrow for a cusp is that the cusp seen as upright by an observer of the plane, and drawn from left to right, receives
an arrow point from left to right. See the remark in section 7 after Theorem 4 for a discussion of this convention. In this coding two consecutive arrows pointing
in the same direction cancel each other, while consecutive arrows with opposing directions do not cancel. State loops with arrow labels can be seen to reduce
uniquely to loops with an irreducible cyclic sequence of arrows. At this point the reader might think that we could use the reduced states to extend even the bracket
for a classical knot diagram (one without any virtual crossings). However, we have the following Proposition.
\bigbreak

\noindent {\bf Proposition 2.} In a classical knot or link diagram, all state loops reduce to
loops that are free from cusps. 
\smallbreak

\noindent {\bf Proof.} The result follows from the Jordan Curve Theorem. Each state loop for a classical knot is a Jordan curve on the surface of the 
$2$-sphere, dividing the surface of the sphere into an inside and an outside. If a state loop has a non-trivial arrow reduction, then there will be non-empty collections of inward-pointing cusps and a collection of outward-pointing cusps
in the reduced loop. Each cusp must be paired with another cusp in the state (since cusps are originally paired, and each reduction can be seen to inherit a pairing
for each cusp as in the reduction procedures of Section 7). No two inward-pointing cusps on a loop can be paired with one another, since they will have incompatible
orientations (similarly for the outward-pointing cusps). Therefore, for a given non-trivial reduced loop $\lambda$ in a classical diagram there must be another non-trivially
reduced loop $\lambda'$ inside  $\lambda$ and another such loop $\lambda''$ outside $\lambda$ to handle the necessary pairings. This means that the given state
of the diagram would have infinitely many loops. Since we work with knot and link diagrams with finitely many crossings, this is not possible. Hence there are no
non-trivially reduced loops in the states of a classical diagram. This completes the proof. //
\bigbreak

\noindent Formally, we have the following state summation for the simple extended bracket $${\cal A}[K] = \Sigma_{S}<K|S>d^{||S||-1} {\cal V}[S]$$
where $S$ runs over the oriented bracket states of the diagram, $<K|S>$ is the usual product of vertex weights as in the 
standard bracket polynomial, $||S||$ is the number of circle graphs in the state $S$, and ${\cal V}[S]$ is a product of the variables $K_{n}$ associated
with the non-trivial circle graphs in the state $S.$ Note that each circle graph (trivial or not) contributes to the power of $d$ in the state summation,
but only non-trivial circle graphs contribute to ${\cal V}[S].$ The regular isotopy invariance of ${\cal A}[K]$ follows from the same analysis that we used for the 
extended bracket, and is combinatorially easier since the reduction rule is simpler.
\bigbreak

\noindent {\bf Theorem 6.} With the above conventions, the simplified extended bracket ${\cal A}[K]$ is a polynomial in $A, A^{-1}$ and the graphical
variables $K_{n}$ (of which finitely many will appear for 
any given virtual knot or link). ${\cal A}[K]$ is a regular isotopy invariant of virtual knots and links. The normalized version
$${\cal W}[K] = (-A^{3})^{-wr(K)} {\cal A}[K]$$ is an invariant virtual isotopy. If we set $A = 1$ and $d = -A^2 - A^{-2} = -2$, then the resulting specialization
$${\cal F}[K] = {\cal A}[K](A = 1)$$
is an invariant of flat virtual knots and links.
\bigbreak

\begin{figure}
     \begin{center}
     \begin{tabular}{c}
     \includegraphics[width=6cm]{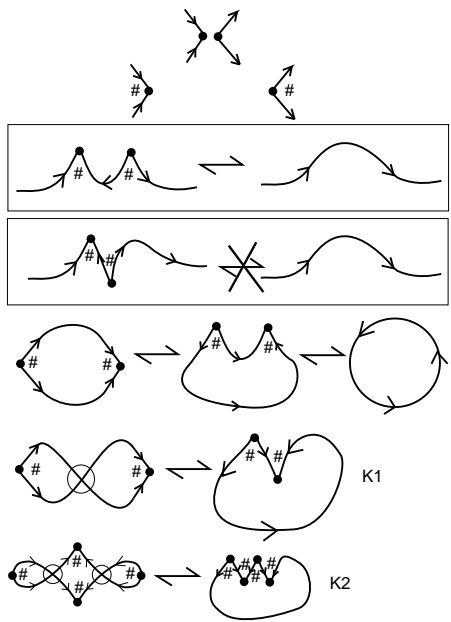}
     \end{tabular}
     \caption{\bf Reduction Relation for Simple Extended Bracket.}
     \label{Figure 48}
\end{center}
\end{figure}

\begin{figure}
     \begin{center}
     \begin{tabular}{c}
     \includegraphics[width=7cm]{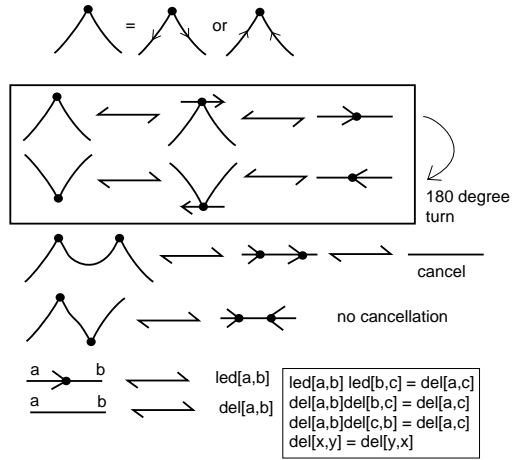}
     \end{tabular}
     \caption{\bf The Arrow Convention.}
     \label{Figure 49}
\end{center}
\end{figure}

\noindent {\bf Example.} Refer to Figures $36$ to $41.$ These depict the extended bracket calculations for the Kishino diagram. It follows at once from
Figure $41$ that if $K$ is the Kishino diagram, then (with $d= -A^2 - A^{-2}$) $${\cal A}[K] = 1 + A^4 + A^{-4} -d^2 K_{1}^{2} + 2 K_{2}.$$ 
It is also easy to see this result directly from the states shown in Figure $36,$ 
since the combinatorics of state reduction is quick for the simple extended bracket. Thus the simple extended bracket shows that the Kishino is 
non-trivial and non-classical. In fact, note that $${\cal F}[K] = 3 + 2 K_{2} - 4 K_{1}^{2}.$$ Thus the invariant ${\cal F}[K]$ of flat virtual diagrams proves that the 
flat Kishino diagram is non-trivial. This example shows the power of the simple extended bracket. 
\bigbreak

\noindent {\bf Remark.} At this writing we do not have an example of a virtual knot $K$ with normalized arrow polynomial ${\cal W}[K]= 1,$ but having a non-trivial extended bracket polynomial.
We conjecture that such examples exist. In the following subsection, we describe a simple program for finding ${\cal A}[K].$
\bigbreak

\subsection{A Mathematica Program for the Arrow Polynomial}

In the lines below we give a concise program for computing the arrow polynomial ${\cal A}[K].$ These program lines end with a sample calculation
for the Kishino diagram. The program is written in 
Mathematica code slightly modified for ease of reading. In Mathematica a function or a replacement rule marks the variables involved in the
replacement. Here we have not marked these variables and so violate the Mathematica syntax. The author of the paper will be happy to supply
working Mathematica code to anyone who is interested in this program. After the program listing, we explain the conventions that drive this program.
These conventions can be viewed in Figures $49$, $50$ and $51.$
\bigbreak
\begin{center}
$$rule1 = \{X[a, b, c, d] :> A del[d, a] del[c, b] + B led[a, b] led[c, d], $$
    $$Y[a, b, c, d] :> B del[d, a] del[c, b] + A led[a, b] led[c, d] \};$$ 
$$rule2 = \{ del[a, b] del[b, c] :> del[a, c],$$ 
$$del[a, b] del[c, b] :> del[a, c], $$
$$del[a, b] led[b, c] :> led[a, c], $$
    $$del[a, b] led[c, b] :> led[c, a], $$
$$del[b, a] led[b, c] :> led[a, c], $$
$$del[b, a] led[c, b] :> led[c, a], $$
    $$led[a, b] del[b, c] :> led[a, c], $$
$$led[a, b] del[c, b] :> led[a, c], $$
$$led[b, a] del[c, b] :> led[a, c], $$
    $$led[b, a] del[b, c] :> led[a, c],$$ 
$$led[a, b] led[b, c] :> del[a, c] \};$$ 
$$rule3 = \{ del[a, a] :> J, $$
$$del[a, b]^2 :> J, $$
$$led[a, b]^2 :> J K1,$$ 
$$led[a, b] led[c, b] led[c, d] led[a, d] :> J K2, $$
    $$led[a, g] led[a, i] led[c, g] led[c, k] led[e, i] led[e, k] :> J K3 \}; $$
$$RawBracket[t] := Simplify[Expand[t /. rule1] //. rule2 /. rule3]$$
$$rule4 = \{ B :> 1/A, J :> -A^2 - 1/A^2 \}; $$
$$Arrow[t] := Simplify[RawBracket[RawBracket[t]/J] /. rule4]$$
$$Flat[t] := Arrow[t] /. A :> 1$$
{\bf Program List}
$$Kishino = X[a,c,h,b]Y[h,b,g,a]X[d,f,c,e]Y[e,g,d,f];$$
$$Expand[Arrow[Kishino]]$$
$$1 + A^{-4} + A^{4} - 2 K1^{2} - K1^{2} A^{-4} - A^{4} K1^{2} + 2 K2$$
{\bf Sample Computation}
\end{center}
\bigbreak

\begin{figure}
     \begin{center}
     \begin{tabular}{c}
     \includegraphics[width=8cm]{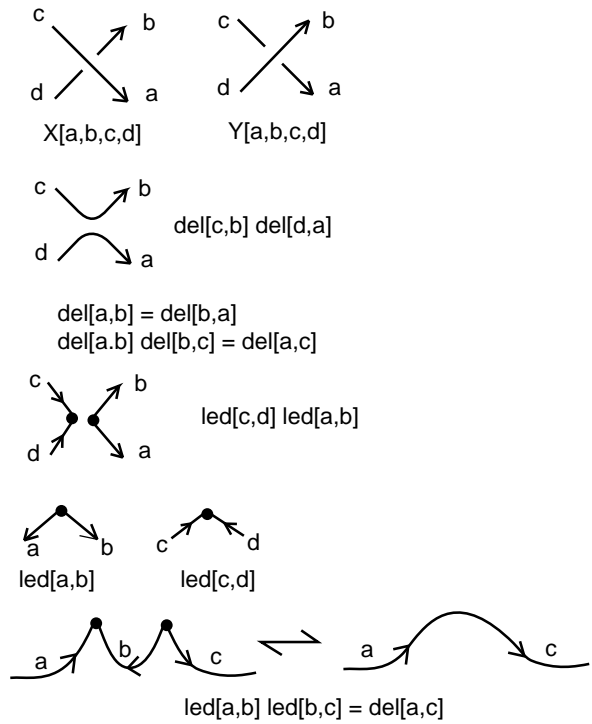}
     \end{tabular}
     \caption{\bf Formal Kronecker Deltas.}
     \label{Figure 50}
\end{center}
\end{figure}

\begin{figure}
     \begin{center}
     \begin{tabular}{c}
     \includegraphics[width=8cm]{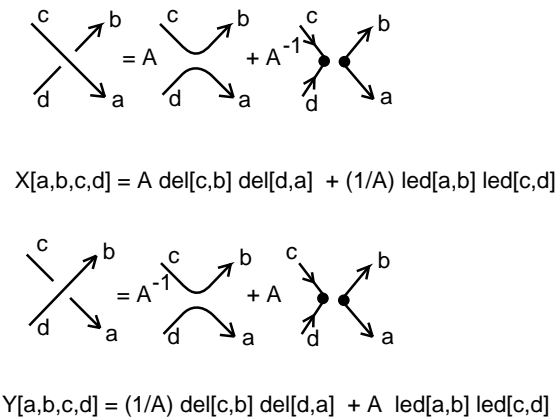}
     \end{tabular}
     \caption{\bf Simple Extended Bracket Expansion via Formal Kronecker Deltas.}
     \label{Figure 51}
\end{center}
\end{figure}

\noindent Figure $50$ illustrates how we take the symbols $X[a,b,c,d]$ and $Y[a,b,c,d]$ as indicators for postive and 
negative crossings respectively. In the underlying $4$-regular graph for the virtual knot diagram each edge from real vertex to real vertex is
given a distinct label. Virtual crossings are not nodes of this graph. Thus each crossing has four labels on its local edges. We take these labels in 
counterclockwise order, starting at the lower right outward pointing arrow when the crossing itself points from left to right. The same convention 
applies to negative crossings. The knot is coded as a formal product of these crossings. Thus in the code lines above we see
$$Kishino = X[a,c,h,b]Y[h,b,g,a]X[d,f,c,e]Y[e,g,d,f].$$
The Kishino diagram is described in terms of its four crossings and how they are connected to one another.
The program formally expands the crossings according to the oriented bracket expansion. Thus we have (see Figure $51$)
$$X[a,b,c,d] \longrightarrow A del[d,a] del[c,b] + B led[a,b] led[c,d]$$ and 
$$Y[a,b,c,d] \longrightarrow  B del[d,a] del[c,b] + A led[a,b] led[c,d].$$
Here $$del[a,b]$$ and $$led[a,b]$$ denote the connections of these labels that occur after the crossing is smoothed in its oriented ($del$) and 
disoriented ($led$) ways. Each of $del[a,b]$ and $led[a,b]$ behave as formal Kronecker deltas, but with different rules. The $del[a,b]$ does not
depend on order, so we have $$del[a,b] del[b,c] = del[a,c]$$ and $$del[a,b] del[c,b] = del[a,c].$$ The $led[a,b]$ does depend upon the order 
{\it by the conventions indicated in Figure $50$}  and we 
only have the ordered reduction $$led[a,b] led[b,c] = del[a,c].$$ This is the symbolic transcription of the basic reduction rule for the simple extended
bracket. See Figure 50 for the comparison with the reduction rules. Note how the two $led$'s combine to form a standard $del.$ Rules $2$ and $3$ in the program
contain all the reductions for $del$ and $led.$ In Rule $3$ we reduce trivial circles to $J$ (which will later become $d = - A^2 - A^{-2}$ and we replace the 
first two non-trivial circle graphs by $K1$ and $K2$ respectively. (Note that in the replacement rules, the circle graphs are
replaced by $J K1$ and $JK2$ respectively. This is because each circle graph in the state expansion contributes the loop value $J$ as well as its
combinatorial value. One obtains the standard bracket polynomial by setting all the variables $K_{n}$ equal to unity.) If one applies the program to larger knots and links, then other non-trivial circle 
graphs will appear, and the program can be modified to name them as well. Rule $4$ replaces $B$ by $A^{-1}$ and $J$ by $- A^2 - A^{-2}.$ 
This program is easily modified to produce the writhe-normalized version of the invariant and to specialize to obtain the flat invariant ${\cal F}[K].$
\bigbreak

\section{Estimating Genus using the Extended Bracket and Arrow polynomials} 
In this section we point out that the genus of the graphs that survive into the expansion of the extended bracket polynomial $<<K>>$ give lower bounds on the 
genus of the virtual link $K.$ Let the genus be denoted by $g(K),$ the least genus surface on which the virtual knot $K$ can be realized. We have the following
theorem. 
\bigbreak

\noindent {\bf Theorem.} Let $G$ be one of the virtual graphs obtained as a reduced state with non-zero coefficient in the expression for the extended bracket polynomial
of a virtual link $K.$ Let $g(G)$ denote the least genus surface supporting the graph $G.$ Then $g(G) \le g(K).$
\bigbreak

\noindent {\bf Proof.} Since a representation of $K$ on a surface will entail the representation (embedding) of each graph $G$ resulting from state
reduction, it follows at once that the genus of each graph $G$ is a lower bound for the genus of $K.$ This completes the proof of the Theorem. //
\bigbreak

\noindent {\bf Remark on Graph Genus.} We can calculate $g(G)$ for any virtual graph $G$ that is given as an immersed graph in the plane
with virtual crossings. The planar data allows the calculation of the genus of the graph, as we shall explain below. This is in marked contrast to the situation
for virtual knots and links, where the flexibility in changing their diagrams by virtual isotopy makes determination of the genus more difficult. 
View Figure 52. In that figure we illustrate a virtual graph and the construction of its embedding in an orientable surface of minimal genus for this graph.
The method indicated in this figure gives an algorithm that can be applies to any virtual graph to determine its genus. One forms the abstract ribbon diagram
$N(G)$ (see the Figure 52) associated to the graph $G$ by replacing each graphical node by a surface neighborhood of that node, and replacing each
virtual crossing by two strips of ribbon, one for each of the arcs running through the virtual crossing. 	These local neighborhoods and ribbon strips are then 
connected by ribbons that form neighborhoods of the remaining arcs in the diagram. Note that from the virtual crossings there are now, in the diagram for
$N(G),$ strips that over and under cross one another. The surface we are actually constructing is an abstract surface, so it does not matter how we represent
these under and over crossing strips. The construction of $N(G)$ as an abstract surface with boundary is unique as a function of the virtual isotopy class of the 
virtual graph $G.$ Now form a surface $S(G)$ by attaching one disk for each boundary component of $N(G).$ We then have an embedding of $G$ in the orientable
surface $S(G),$ and it follows easily that $S(G)$ is a surface of least genus in which the virtual graph can be embedded. The reason for this minimality is simply
that any surface in which $G$ is embedded will contain a neighborhood of $G$ that is homeomorphic to $N(G).$ Thus adding disks to each boundary component of 
$S(G)$ yields the least genus surface. Furthermore, we can calculate the genus of $S(G)$ in terms of the graph $G$ as follows: Let $v(G)$ denote the number of
vertices of $G$ and $e(G)$ denote the number of edges of $G.$  Let $f(G)$ denote the number of boundary components of $N(G).$ Then we have, Euler's formula
$$2 - 2g(S(G)) = v(G) - e(G) + f(G).$$ Note also, that since $G$ is locally of degree $4$, we have the $e(G) = 2v(G).$ 
From this we obtain the formula for the genus of $S(G).$ $$g(G) = g(S(G)) = 1 + \frac{v(G) - f(G)}{2}.$$ For example, in Figure 52 we have
$v(G) = 2 = f(G)$ and so $g(G) = 1.$
\bigbreak

\noindent {\bf Examples.} In Figure 53 we illustrate a family of flat virtual links $L_{n}$ for $n$ a natural number. It is not hard to see
that the virtual graph $G_{n}$ obtained from $L_{n}$ by making each virtual crossing into a graphical node survives in the extended bracket
invariant of $L_{n}.$ This means that $K_{1}^{2} K_{2}^{n-1}$ survives in the arrow polynomial. Furthermore, a calculation with the graphs
$G_{n}$, using the technique of the remark above, shows that $g(G_{n}) = n.$ Thus, the Theorem on genus tells us that $g(L_{n}) = n$
for $n=1,2,\cdots n.$ See Figure 54 for an illustration of the combinatorics supporting these statements for $n = 2.$ In Figure 55 we illustrate
a flat virtual knot $K$ and a state $S$ of the extended bracket expansion of $K$. The state $S$ reduces to the state $T$ with graph $G$ as shown
in this figure. By using the method of the remark above, it follows that $G$ has genus two. Since it is easy to see that $K$ can be represented in 
a surface of genus two, it follows that $g(K) = 2.$ Note that in this last example, the genus estimate comes from the reduced state and that this state does 
survive in the extended bracket invariant. \bigbreak

\begin{figure}
     \begin{center}
     \begin{tabular}{c}
     \includegraphics[width=7cm]{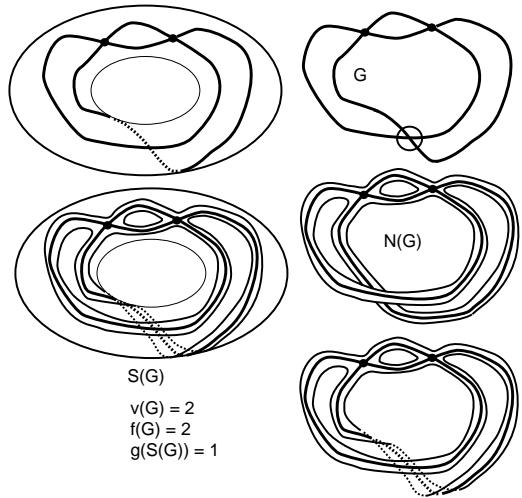}
     \end{tabular}
     \caption{\bf Genus of Virtual Graphs.}
     \label{Figure 52}
\end{center}
\end{figure}

\begin{figure}
     \begin{center}
     \begin{tabular}{c}
     \includegraphics[width=7cm]{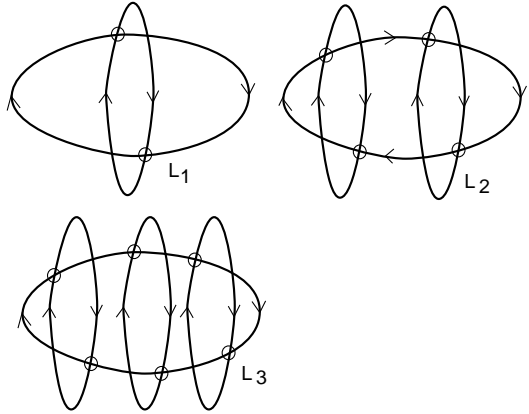}
     \end{tabular}
     \caption{\bf The Family $L_{n}$.}
     \label{Figure 53}
\end{center}
\end{figure}

\begin{figure}
     \begin{center}
     \begin{tabular}{c}
     \includegraphics[width=7cm]{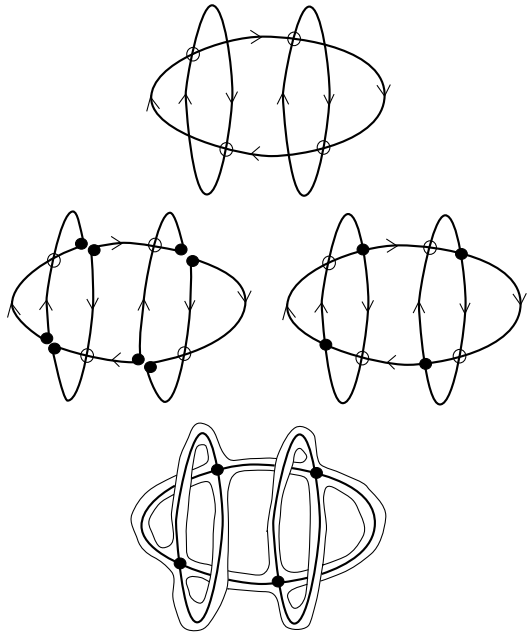}
     \end{tabular}
     \caption{\bf The Flat Link $L_{2}$.}
     \label{Figure 54}
\end{center}
\end{figure}

\begin{figure}
     \begin{center}
     \begin{tabular}{c}
     \includegraphics[width=7cm]{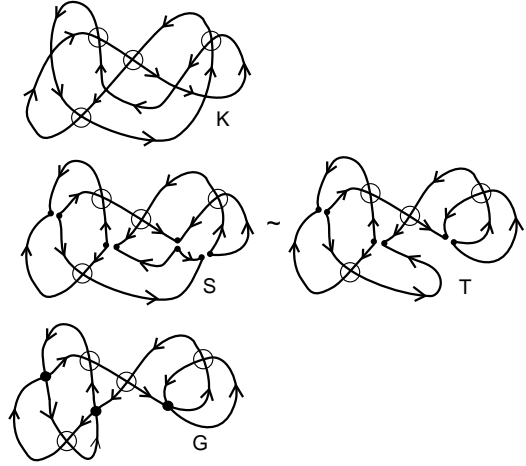}
     \end{tabular}
     \caption{\bf A Flat Knot of Genus Two.}
     \label{Figure 55}
\end{center}
\end{figure}

\section {Discussion} 
It is a key problem to understand the full
extent to which the extended bracket polynomial can detect non-classicality for non-trivial virtual knots with unit Jones polynomial. Recall that the
extended bracket is an invariant for long virtual knots and (by taking $A = 1$ and $d = -2$) for flat virtual links and long flat virtual
knots. It is a feature of the method we use
to  associate virtual graphs to states, that the extended bracket can detect some long flat links whose closures are trivial. Thus the extended 
bracket can be used to investigate the kernel of the closure mapping from long virtuals to closed viruals. In this paper we have given a number of
specific diagrammatic calculations of the extended bracket invariant. More systematic methods of calculation for the extended bracket are needed.
We have also given a simplified version, ${\cal A}[K],$ of the extended bracket that is a state sum with an infinite set of extra polynomial variables.
The simple extended bracket is a strong invariant of virtual links and flat virtual diagrams and is easily computable by a Mathematica program  described 
in this paper.
\bigbreak

A number of natural variants of the extended bracket polynomial are available. For example, we can obtain a powerful invariant of long flats
$F$ by computing \,\,\,\,\,\, $<<A(F)>>,$ using the Embedding Theorem that we described in the section on long virtuals. Furthermore, given a knot in a
specific thickened surface $\Sigma$, we modify the definition of $<<K>>$ so that the curves associated with the states are homotopy classes of
immersions in the surface $\Sigma.$ The same idea applies to a curve that is immersed in $\Sigma,$ representing a flat virtual knot. There may be further
variants. For example, it is possible that one can apply a similar procedure to extend the Miyazawa \cite{Miyazawa1,Miyazawa2,Kamada1,Kamada2,KamadaMiya}
polynomials. It is natural to ask for a categorification \cite{Kho,DB1,DB2,M,V} of the extended bracket invariant. 
\bigbreak

\end{document}